\documentclass[reqno]{amsart}

\usepackage{amssymb,young}
\usepackage{array, boldline, makecell, booktabs}
\usepackage{bm}
\usepackage{ytableau}
\usepackage{euscript, mathrsfs}
\usepackage{tikz}
\usetikzlibrary{backgrounds}
\usepackage{mathtools}
\usepackage{cite}
\usepackage[colorlinks]{hyperref}
\usepackage{cleveref}
\usepackage[
    margin=1in,
    headsep=14pt,
    footskip=0.3in
]{geometry}
\numberwithin{equation}{section}    
\usepackage[all]{xy}
\usepackage{graphicx,import}
\usepackage{amscd}
\usepackage{color}
\usepackage{enumerate}
\usepackage{fourier}
\usepackage{cases}
\usepackage[T1]{fontenc}

\usetikzlibrary{arrows.meta}

\newtheorem{mainthm}{Theorem}

\newtheorem{thm}{Theorem}[section]
\newtheorem{lem}[thm]{Lemma}
\newtheorem{proposition}[thm]{Proposition}
\newtheorem{corollary}[thm]{Corollary}
\theoremstyle{definition}
\newtheorem{example}[thm]{Example}
\newtheorem{definition}[thm]{Definition}
\newtheorem{conj}[thm]{Conjecture}
\newtheorem{rmk}[thm]{Remark}

\DeclareMathOperator{\aut}{aut}

\DeclareMathOperator{\cPF}{cPF}

\DeclareMathOperator{\pdinv}{pdinv}
\DeclareMathOperator{\ldinv}{ldinv}
\DeclareMathOperator{\cont}{cont}

\DeclareMathOperator{\Hessen}{Hess}
\DeclareMathOperator{\aseq}{aseq}

\DeclareMathOperator{\h}{\mathfrak{h}}

\DeclareMathOperator{\pib}{\pi^{\bullet}}
\DeclareMathOperator{\pLD}{pLD}
\DeclareMathOperator{\cLLT}{cLLT}
\DeclareMathOperator{\z}{\mathbf{z}}
\DeclareMathOperator{\mm}{\mathbf{m}}
\DeclareMathOperator{\sLD}{sLD}
\DeclareMathOperator{\LD}{LD}
\DeclareMathOperator{\skel}{Skel}
\DeclareMathOperator{\asc}{asc}

\DeclareMathOperator{\CH}{\mathbf{H}}
\DeclareMathOperator{\CE}{\mathbf{E}}

\DeclareMathOperator{\area}{area}

\DeclareMathOperator{\std}{std}

\DeclareMathOperator{\inv}{inv}
\DeclareMathOperator{\dinv}{dinv}

\DeclareMathOperator{\Set}{Set}

\DeclareMathOperator{\stat}{stat}

\DeclareMathOperator{\pos}{pos}

\newcommand{\mmb}{\mathbf{m}^{\bullet}}

\newcommand{\Lam}{\Lambda}

\newcommand{\Eqt}{\mathcal{E}_{q,t}}

\title{The elliptic Hall algebra and the combinatorics of affine Springer fibers}

\author{Donghyun Kim}
\address{Department of Mathematics \\ Ewha Womans University \\Seoul  \\ Korea}
\email{kdh310@ewha.ac.kr}

\author{Jaeseong Oh}
\address{Department of Mathematics \\ Sungkyunkwan University \\ Suwon \\ Korea}
\email{jaeseongoh@skku.edu}


\begin{document}

\begin{abstract}
For coprime positive integers $m,n$ and $k\ge1$, we study
the symmetric functions $e_k[-X^{m,n}Y]\cdot1$ arising
from the elliptic Hall algebra action.
We identify their $\omega_X$-transforms with the normalized
bigraded Frobenius characteristics computed by
Boixeda Alvarez--Hameister--Ma for the torus-equivariant
Borel--Moore homology of the positive parts of homogeneous
affine Springer fibers.
We introduce cyclic LLT polynomials and determine their
monomial expansions, yielding explicit monomial formulas
for these elliptic Hall algebra expressions in terms of
cyclic parking functions.
These results extend the affine Springer fiber identity
and monomial formula of Carlsson--Mellit from integral
slopes to arbitrary positive rational slopes.
\end{abstract}

\maketitle

\section{Introduction}\label{Sec: Intro}

The elliptic Hall algebra, introduced by Burban and
Schiffmann~\cite{BS12}, has played a central role in the study
of shuffle theorems, diagonal coinvariant spaces, and the
geometry of Hilbert schemes and affine Springer fibers.
It contains distinguished copies $\Lambda(X^{m,n})$ of the
ring $\Lambda$ of symmetric functions, indexed by coprime
pairs $(m,n)$, and acts on $\Lambda$.
Schiffmann--Vasserot realize this action by correspondences on
the localized $(\mathbb C^\times)^2$-equivariant $K$-theory of
Hilbert schemes of points on $\mathbb C^2$~\cite{SV13}.
This action also underlies the rational shuffle theorem
and its variants, which give combinatorial formulas for
symmetric functions of the form $f[-MX^{m,n}]\cdot1$, where
$f\in\Lambda$ and
$M=(1-q)(1-t)$~\cite{CM18, Mel21, BHMPS23, BHMPS21Delta, BHMPS25LW}.

For coprime positive integers $m,n$ and $k\ge1$, we study
the symmetric functions
\begin{equation}\label{eq: EHA operation}
    e_k\!\left[-X^{m,n}Y\right]\cdot1
    =
    \sum_{\mu\vdash k}
    \left(e_\mu\!\left[-X^{m,n}\right]\cdot1\right)m_\mu[Y],
\end{equation}
where the expansion follows from the dual Cauchy identity.
When $n=1$, this recovers the symmetric function appearing in
Carlsson--Mellit's work on the nabla operator~\cite{CM25}
and affine Springer fibers~\cite{CM21}.
In a companion paper, Boixeda Alvarez--Hameister--Ma compute
equivariant $K$-theory classes of coherent sheaves on Hilbert
schemes in terms of elliptic Hall algebra
operators~\cite[Theorem~9.5]{BAHM26}.

In this paper, we relate the symmetric functions
in~\eqref{eq: EHA operation}
to the geometry of affine Springer fibers. Let $\mathcal B_{km,kn}$ be the positive part of a regular
semisimple homogeneous affine Springer fiber for $GL_{kn}$
of slope $m/n$.
Its Borel--Moore homology carries commuting $S_{kn}$- and
$S_k$-actions, arising from Springer theory and monodromy,
respectively.
We use the alphabets $X$ and $Y$ to record these two actions.
Boixeda Alvarez--Hameister--Ma compute the joint Frobenius
characteristic of the centralizer-torus-equivariant
Borel--Moore homology of $\mathcal B_{km,kn}$
\cite[Corollary~5.25]{BAHM26}.
In our parametrization and grading conventions, the normalized
equivariant Frobenius characteristic takes the form
\begin{equation}\label{eq: affine Springer fiber formula}
    \mathcal F_{km,kn}[X,Y;q,t]
    \coloneqq
    \sum_{\lambda\vdash kn}
    m_\lambda[X]
    \sum_{\substack{
        \pi^\bullet\in\overline{\cPF}_{m,n}^k(\lambda)\\
        \pi^\bullet\text{ dominant}
    }}
    t^{\area(\pi^\bullet)}
    q^{\overline{\stat}(\pi^\bullet)}
    \frac{
        \Xi_{\mmb(\pi^\bullet)}
        \!\left[\dfrac{Y}{1-q};q\right]
    }{
        [\alpha'(\pi^\bullet)]_q!
    }.
\end{equation}
Here the inner sum runs over dominant tuples of weak cyclic
parking functions, and $\Xi_{\mmb(\pi^\bullet)}$ is the associated
cyclic chromatic symmetric function.
See Section~\ref{sec:affine-springer} for the grading normalization
and Sections~\ref{Sec: cPF} and~\ref{Sec: asf combinatorics}
for the combinatorial definitions.

Our first main theorem, combined with the geometric character
formula in~\cite[Corollary~5.25]{BAHM26}, extends the Frobenius
characteristic formula of Carlsson--Mellit~\cite[Theorem~A(b)]{CM21}
from the $n=1$ case to arbitrary coprime positive
integers $m,n$.

\begin{mainthm}\label{thm: main}
For coprime positive integers $m,n$ and every $k\ge1$,
\[
    \omega_X\!\left(e_k[-X^{m,n}Y]\cdot 1\right)
    =
    \mathcal F_{km,kn}[X,Y;q,t].
\]
\end{mainthm}

Theorem~\ref{thm: main} supplies the combinatorial bridge
needed to complete the comparison in the companion
paper~\cite[Theorems~1.1 and~9.7]{BAHM26}.
This comparison identifies the normalized bigraded characters
of the ordinary Borel--Moore homology of $\mathcal B_{km,kn}$
and of the associated graded of rational Cherednik algebra
modules of minimal support with the corresponding equivariant
Euler characteristics on centered Hilbert schemes.
Together with the torus-link shuffle theorem~\cite{KO26},
these results establish the torus-link extension of the
Gorsky--Oblomkov--Rasmussen--Shende
correspondence~\cite{GORS14}.

The formula for $\mathcal F_{km,kn}$ is a monomial expansion
in $X$, but its coefficients involve cyclic chromatic symmetric
functions evaluated at the plethystic alphabet $Y/(1-q)$.
To obtain an explicit monomial formula in both alphabets,
we therefore study the monomial expansions of these
plethystic transforms.

For an $n$-tuple $\mmb$ of Hessenberg functions of size $k$,
let $\Xi_{\mmb}[Y;q]$ be the associated cyclic chromatic
symmetric function. We define the cyclic LLT polynomial by
\[
    \cLLT_{\mmb}[Y;q]
    \coloneqq
    (1-q)^k\,
    \Xi_{\mmb}\!\left[\frac{Y}{1-q};q\right].
\]
When $n=1$, this recovers the usual unicellular LLT polynomial
associated with a Hessenberg function.
In Theorem~\ref{thm: rLLt monomial}, we express its monomial
coefficients as weighted sums over sorted loop decompositions.
We also establish the corresponding expansion for the
plethystic transform of the dual cyclic chromatic symmetric
function.
Together with Theorem~\ref{thm: main} and the involution
$\omega_X$, these expansions yield the following monomial
formula for $e_k[-X^{m,n}Y]\cdot1$, generalizing
\cite[Theorem~A]{CM25}.

    \begin{mainthm}[Theorem~\ref{thm: cm general}]
\label{thm: main B} For coprime positive integers $m,n$ and every $k\ge1$, we have
\begin{align}\label{eq: combinatorial formula}
    e_k\left[-X^{m,n} Y\right]\cdot1
    &=
    \sum_{\substack{
        \left(\pi^\bullet,b\right)\in
        \cPF_{m,n}^k\times\mathbb Z_{\ge1}^k\\
        \left(\pi^\bullet,b\right)\text{ sorted}
    }}
    \frac{
        t^{\area\left(\pi^\bullet\right)}
        q^{\stat\left(\pi^\bullet\right)}
    }{
        \left(1-q\right)^k\aut_q\left(\pi^\bullet,b\right)
    }
    X_{\pi^\bullet}Y_b.
\end{align}
\end{mainthm}
Here $X_{\pi^\bullet}$ and $Y_b$ record the labels of the
parking-function tuple and the auxiliary labels $b$,
respectively.
The sorting condition, statistics, and factor $\aut_q$ are
defined in Sections~\ref{Sec: cPF}
and~\ref{Sec: rational uLLT}.

The paper is organized as follows.
Section~\ref{Sec: Pre} recalls the necessary notation for
symmetric functions, plethystic substitution, the elliptic Hall
algebra, and affine Springer fibers.
Section~\ref{Sec: cPF} defines cyclic parking functions and their statistics.
Section~\ref{Sec: h operators} develops the
$\mathfrak h$-operators and their weak analogues,
introduces the weak $\mathfrak e$-operators, and derives
the operator expansion used in Theorem~\ref{thm: main}.
Section~\ref{Sec: asf combinatorics} compares this expansion
with the affine Springer fiber formula via cyclic chromatic
symmetric functions and completes the proof of
Theorem~\ref{thm: main}.
Finally, Section~\ref{Sec: rational uLLT} introduces cyclic LLT
polynomials, proves their monomial expansions, and establishes
Theorem~\ref{thm: main B}. Appendix~\ref{Sec: Appendix} compares cyclic parking functions
with affine permutations and gives a Hilbert-series 
formula.

\subsection*{Acknowledgments}
The authors thank Pablo Boixeda Alvarez, Thomas Hameister,
and Xinchun Ma for sharing their formula and explaining
the underlying geometric picture. They also thank
Anton Mellit for a fruitful discussion.

\subsection*{Declaration on Generative AI}
The authors used the ChatGPT 5.5, 5.6, and 6 Pro models
to improve the language and exposition of the manuscript.
All mathematical content is solely the work of the authors.

\section{Preliminaries}\label{Sec: Pre}
For an integer $r$, write $[r]=\{i\in\mathbb Z:1\leq i\leq r\}$;
in particular, $[r]=\varnothing$ for $r\leq0$.
We use the standard notation for $q$-analogues:
\[
[n]_q=\sum_{i=0}^{n-1}q^i,\qquad
[n]_q!=\prod_{i=1}^n[i]_q\quad(n\geq1),\qquad
[0]_q!=1.
\]
\subsection{Symmetric functions and the elliptic Hall algebra}
\label{subsec:symf-eha}
 A \emph{composition} of an integer $n\geq 0$ is a finite
sequence $\alpha=\left(\alpha_1,\ldots,\alpha_\ell\right)$ of positive integers such that
$\left|\alpha\right|\coloneqq\alpha_1+\cdots+\alpha_\ell=n$; for $n=0$, we allow the empty composition. We write \(\alpha\models n\) and call \(\ell(\alpha)\coloneqq\ell\) its length. The entries $\alpha_i$ are referred to as the \emph{parts} of
$\alpha$. We define $\Set\left(\alpha\right)
\coloneqq
\left\{
    \sum_{i=1}^{j}\alpha_i:
    1\le j<\ell(\alpha)
\right\}$ and $\Set\left(\alpha\right)^c
    =
    \left[n-1\right]\setminus \Set\left(\alpha\right)$. We write $\left[\alpha\right]_{q}!\coloneqq\prod_{i=1}^{\ell}\left[\alpha_i\right]_q!$. For example, if $\alpha=\left(2,3,1\right)\models 6$, then $\Set\left(\alpha\right)=\left\{2,5\right\}$ and $\left[\alpha\right]_{q}!=\left[2\right]_{q}!\left[3\right]_{q}!\left[1\right]_{q}!=\left(1+q\right)^2 \left(1+q+q^2\right)$. For two compositions $\alpha,\beta$, we write $\alpha \sim \beta$ if one is a rearrangement of the other.

A \emph{partition} of $n$, denoted by $\lambda\vdash n$, is a composition
$\lambda=\left(\lambda_1,\ldots,\lambda_\ell\right)$ such that
$\lambda_1\geq\lambda_2\geq\cdots\geq\lambda_\ell>0$.
The \emph{Young diagram} of $\lambda$ is
\[
    \left\{\left(i,j\right)\in\mathbb Z^2:
      1\le i\le \ell\left(\lambda\right),\ 1\le j\le \lambda_i\right\}.
\]
We identify a partition with its Young diagram.

Let $\Lam$ denote the ring of symmetric functions over $\mathbb Q\left(q,t\right)$, and let
$X=x_1+x_2+\cdots$ be the usual infinite alphabet. We follow the usual convention for the bases of $\Lam$ indexed by partitions: the monomial, elementary,
complete homogeneous, power-sum, and Schur symmetric functions, denoted
respectively by $m_\lambda$, $e_\lambda$, $h_\lambda$, $p_\lambda$, and $s_\lambda$.
The involution $\omega$ is determined by $\omega\left(e_n\right)=h_n$, or
equivalently by $\omega\left(p_k\right)=\left(-1\right)^{k-1}p_k$. When applying this involution to symmetric functions in the
$X$-variables, we write it as $\omega_X$.

The \emph{monomial quasisymmetric function} $M_\alpha$ indexed by
a composition $\alpha=\left(\alpha_1,\ldots,\alpha_\ell\right)$ is defined by
\[
    M_\alpha
    \coloneqq
    \sum_{i_1<i_2<\cdots<i_\ell}
    x_{i_1}^{\alpha_1}x_{i_2}^{\alpha_2}\cdots x_{i_\ell}^{\alpha_\ell}.
\]
Then
\[
    m_\lambda=\sum_{\alpha\sim\lambda}M_\alpha.
\]
Thus an expression $f=\sum_\alpha c_\alpha M_\alpha$ is symmetric if and only if
$c_\alpha=c_\beta$ whenever $\alpha\sim\beta$.

The \emph{Hall inner product} on $\Lam$ is defined by
\[
    \left\langle h_\lambda,m_\mu\right\rangle=\delta_{\lambda\mu};
\]
that is, the complete homogeneous and monomial bases are dual to each other.
For $f\in\Lam$, let $f^\perp$ denote the adjoint of multiplication by $f$:
\[
    \left\langle f^\perp g,h\right\rangle=\left\langle g,fh\right\rangle
    \qquad \left(g,h\in\Lam\right).
\]

We employ \emph{plethystic notation} throughout. An \emph{alphabet} is regarded as a formal sum of
variables, monomials, or virtual terms. Plethystic substitution is determined linearly and multiplicatively
by the rules $p_k\left[A+B\right]=p_k\left[A\right]+p_k\left[B\right]$, $p_k\left[AB\right]=p_k\left[A\right]p_k\left[B\right]$, $p_k\left[u\right]=u^k$ for a
monomial $u$, and $p_k\left[A-B\right]=p_k\left[A\right]-p_k\left[B\right]$. Furthermore, we set
$M=\left(1-q\right)\left(1-t\right)$; as an example,
$p_k\left[-X/M\right]=-p_k\left[X\right]/\left(\left(1-q^k\right)\left(1-t^k\right)\right)$.

We now briefly recall the action of the \emph{elliptic Hall algebra} used in the sequel.
Let $\Eqt$ be the elliptic Hall algebra over $\mathbb Q\left(q,t\right)$, with the
normalization from \cite{BS12} and the notation established in \cite{BHMPS21Delta}.
For every coprime pair $\left(m,n\right)\in\mathbb Z^2$, the algebra contains
a copy of $\Lam$, denoted by $\Lam\left(X^{m,n}\right)$; the element corresponding to $f\in\Lam$ is written as $f\left(X^{m,n}\right)$. The algebra also contains the central Laurent
polynomial algebra $\mathbb Q\left(q,t\right)\left[c_1^{\pm1},c_2^{\pm1}\right]$. We omit
the defining relations here, as they are fully detailed in \cite[Section~3.3]{BHMPS21Delta}.
If $E$ is a plethystic expression, $f\left[EX^{m,n}\right]$ denotes the element of
$\Lam\left(X^{m,n}\right)$ that corresponds to $f\left[EX\right]$ under this identification.

The representation of $\Eqt$ on $\Lam$ is explicitly described on the four standard
primitive directions as follows. Let $\left\{\widetilde H_\mu\left[X;q,t\right]\right\}_\mu$ be the
\emph{modified Macdonald} basis for $\Lambda$. The \emph{diagram generator} associated with $\mu$ is defined as
$B_\mu\left(q,t\right)\coloneqq\sum_{\left(i,j\right)\in\mu}q^{j-1}t^{i-1}$. If $A=A\left(q,t\right)$, let
$\overline A\coloneqq A\left(q^{-1},t^{-1}\right)$. We define diagonal operators by
$f\left[B\right]\widetilde H_\mu=f\left[B_\mu\left(q,t\right)\right]\widetilde H_\mu$ and
$f\left[\bar B\right]\widetilde H_\mu=f\left[\overline{B_\mu\left(q,t\right)}\right]\widetilde H_\mu$; additionally,
$g\left[X\right]^\bullet$ denotes the operator of multiplication by $g\left[X\right]$.

\begin{proposition}[{\cite[Proposition~3.3.1]{BHMPS23}}]
\label{prop:eha-action}
Under the above action of $\Eqt$ on $\Lam$, one has
\begin{equation}\label{eq:eha-standard-action}
\begin{aligned}
    c_1&\mapsto 1, &
    c_2&\mapsto \left(qt\right)^{-1}, \\
    f\left(X^{1,0}\right)&\mapsto \left(\omega f\right)\!\left[B-\frac{1}{M}\right], &
    f\left(X^{-1,0}\right)&\mapsto \left(\omega f\right)\!\left[\overline{\frac{1}{M}-B}\right], \\
    f\left(X^{0,1}\right)&\mapsto f\!\left[-\frac{X}{M}\right]^\bullet, &
    f\left(X^{0,-1}\right)&\mapsto f\left[X\right]^\perp .
\end{aligned}
\end{equation}
\end{proposition}
Thus, vertical generators act via multiplication and Hall adjoints, whereas
horizontal generators act diagonally on the modified Macdonald basis.

Finally, we recall the \emph{nabla operator} $\nabla$ introduced in \cite{BGHT99}. It is
defined by the relation
\begin{equation}\label{eq:nabla-definition}
    \nabla\widetilde H_\mu\left[X;q,t\right]
    =T_{\mu}\widetilde H_\mu\left[X;q,t\right],
\end{equation}
where $T_{\mu}=\prod_{\left(i,j\right)\in\mu}q^{j-1}t^{i-1}$.
Although $\nabla$ is not intrinsically an element of $\Eqt$, conjugation by $\nabla$ serves to shift the slope of the
above action.

\begin{lem}[{\cite[Lemma~3.4.1]{BHMPS23}}]
\label{lem:nabla-conjugation}
For every coprime pair $\left(m,n\right)$ and every $f\in\Lam$,
\begin{equation}\label{eq:nabla-conjugation}
    \nabla f\left(X^{m,n}\right)\nabla^{-1}=f\left(X^{m+n,n}\right)
\end{equation}
as operators on $\Lam$.
\end{lem}

Combining Proposition~\ref{prop:eha-action} and Lemma~\ref{lem:nabla-conjugation}, we obtain
\begin{equation}\label{eq:eha-nabla-identity}
    f\left[-M X^{m,1}\right]\cdot 1=\nabla^m f\left[X\right]
    \qquad \left(m\in\mathbb Z\right).
\end{equation}
The case $m=0$ is exactly the vertical action, and the general case follows via
conjugation by $\nabla^m$. Consequently, the symmetric functions that appear in various types of
shuffle theorems can be interpreted as elliptic Hall algebra
operators applied to the vacuum vector $1$. More generally, using the identification in \cite[Proposition~6.7]{Neg14}, we can interpret the rational shuffle theorem \cite{CM18, Mel21}, the rational Loehr--Warrington-type theorem
\cite{KO24, BHMPS25LW}, and related constructions as combinatorial formulas for the symmetric functions
\[
    e_k\left[-M X^{m,n}\right]\cdot 1, \qquad 
    s_\lambda\left[-M X^{m,n}\right]\cdot 1
\]
for suitably chosen coprime pairs $\left(m,n\right)$ and relevant indices.

\subsection{Affine Springer fibers}
\label{sec:affine-springer}

Let $m,n$ be coprime positive integers and let $k\ge1$.
Put $K=\mathbb C((z))$ and $\mathcal O=\mathbb C[[z]]$.
For $G=GL_{kn}$ and an Iwahori subgroup
$I\subset G(\mathcal O)$, the affine Springer fiber
associated with $\gamma\in\mathfrak{gl}_{kn}(K)$ is
\[
  \mathrm{Fl}_\gamma
  =
  \left\{
    gI\in G(K)/I:
    \operatorname{Ad}(g^{-1})\gamma
    \in\operatorname{Lie}(I)
  \right\}.
\]
We take $\gamma$ to be regular semisimple and homogeneous
of slope $m/n$. Here homogeneity means that
$s^{-km}\gamma(s^{kn}z)$ is conjugate to $\gamma(z)$
for every $s\in\mathbb C^\times$.

We write $\mathcal B_{km,kn}=\mathrm{Fl}_\gamma^+$ for its
positive part, following \cite[Definition~5.9]{BAHM26}.
This is the closed union of strata corresponding to
positive $(km,kn)$-invariant sets: cofinite subsets
$\Delta\subset\mathbb Z_{>0}$ satisfying
$\Delta+km\subset\Delta$ and $\Delta+kn\subset\Delta$.
The statistic
\[
  \operatorname{area}(\Delta)
  =
  \#(\mathbb Z_{>0}\setminus\Delta)
\]
defines the Hikita filtration $F_\bullet$ on the
Borel--Moore homology of $\mathcal B_{km,kn}$.
The closed area truncations are finite-dimensional
projective varieties, and the characters below are
understood using these truncations.

There are commuting actions of $S_{kn}$ and $S_k$ on this
homology. The first is the Springer action; the second
arises from monodromy in the family of regular semisimple
homogeneous elements of the given slope.
We use the alphabets $X$ and $Y$ to record their respective
Frobenius characteristics.
The Hessenberg paving of
Goresky--Kottwitz--MacPherson \cite{GKM06}, developed in
\cite[Section~5]{BAHM26}, expresses the relevant strata as
affine bundles over generalized cyclic Hessenberg
varieties. Their monodromy characters are cyclic chromatic
symmetric functions. Combining these characters with the
affine-bundle degrees gives the ordinary joint character
computed in \cite[Theorem~5.24]{BAHM26}.

Following the grading conventions of~\cite{BAHM26}, write
the ordinary joint character as
\[
  \mathcal G_{km,kn}[X,Y;q,t]
  =
  \sum_{r,i\geq0} q^r t^i\,
  \operatorname{Frob}_{S_{kn}\times S_k}
  \left(
    \operatorname{gr}_r^F
    H^{\mathrm{BM}}_{2i}(\mathcal B_{km,kn})
  \right).
\]
Here $q$ records the Hikita filtration and $t$ records
half the homological degree.

We retain equivariance for the rank-$k$ block-scaling torus
$T_\gamma\cong(\mathbb C^\times)^k$, including its diagonal
scalar subgroup.
The corresponding joint Frobenius characteristic is computed
in~\cite[Corollary~5.25]{BAHM26}.
In these grading conventions, the equivariant refinement
is expressed by the plethystic substitution
$Y\mapsto Y/(1-t^{-1})$; see also
\cite[Corollary~3.30 and equation~(16)]{BAHM26}.
After exchanging the grading variables and reversing the
homological grading, this alphabet becomes $Y/(1-q)$.

Set
\[
  \delta
  =
  \frac{(km-1)(kn-1)+k-1}{2}.
\]
We use the normalization in which ordinary
$H^{\mathrm{BM}}_{2i}$ has degree $q^{\delta-i}$ and the
Hikita filtration has degree $t$.
The corresponding normalized equivariant character is
\begin{align}
  \mathcal F_{km,kn}[X,Y;q,t]
  &=
  \sum_{\substack{r\geq0\\ i\in\mathbb Z}}
  t^r q^{\delta-i}\,
  \operatorname{Frob}_{S_{kn}\times S_k}
  \left(
    \operatorname{gr}_r^F
    H^{T_\gamma,\mathrm{BM}}_{2i}
    (\mathcal B_{km,kn})
  \right)
  \notag\\
  &=
  q^\delta\,
  \mathcal G_{km,kn}
  \left[
    X,\frac{Y}{1-q};t,q^{-1}
  \right].
  \label{eq:equivariant-normalization}
\end{align}
Multiplication by an equivariant parameter lowers
Borel--Moore degree by two and therefore raises
$q$-degree by one. Thus $\mathcal F_{km,kn}$ records the equivariant character
computed in~\cite[Corollary~5.25]{BAHM26},
with the grading normalization
in~\eqref{eq:equivariant-normalization}.

Under the cyclic parking-function parametrization of
Sections~3 and~5, this is the normalization of the series
in~\eqref{eq: affine Springer fiber formula}.
Theorem~\ref{thm: main} identifies this normalized character with
\[
  \omega_X\left(e_k[-X^{m,n}Y]\cdot1\right).
\]

\section{Cyclic parking functions}\label{Sec: cPF}

\subsection{Cyclic parking functions}
\label{subsec:cyclic-parking-functions}

An \emph{$\left(m,n\right)$-path} is a lattice path $P$ consisting of exactly $n$ north
steps and $m$ east steps, starting at $\left(-r,0\right)$ for some
$r\in\mathbb Z_{\geq 0}$, beginning with a north step, and staying weakly above
the line $my=nx$. Thus $P$ ends at $\left(-r+m,n\right)$. The \emph{area sequence} of $P$, denoted by $\aseq\left(P\right)$, is a sequence $\left(a_1,a_2,\dots,a_n\right)$ such that $a_i$ counts the number of full cells to the right of the $i$-th north step and to the left of the line $my=nx$.

Let $N_1,\ldots,N_n$ be the north steps of $P$, listed from bottom to top along
the path. A \emph{labeling} of $P$ is a map $f$ from its north steps to
$\mathbb Z_{>0}$. To state the cyclic condition, we set
\[
    N_{n+1}\coloneqq N_1+\left(m,n\right), \qquad f\left(N_{n+1}\right)\coloneqq f\left(N_1\right).
\]
For the labeled path $\pi=\left(P,f\right)$, write
\[
 \pi_\ell\coloneqq\left(N_\ell,f\left(N_\ell\right)\right)\qquad\left(1\le\ell\le n\right)
\]
for its $\ell$-th \emph{labeled north step}.
We also set $\pi_{n+1}\coloneqq\left(N_{n+1},f\left(N_1\right)\right)$.
When comparing labeled north steps from different sheets,
we compare their positions and labels without retaining the sheet index.
Two labeled north steps with identical positions and labels
are said to have the same labeled north-step type.

A labeling $f$ is called \emph{strictly cyclic} if
$f\left(N_i\right)<f\left(N_{i+1}\right)$, and \emph{weakly cyclic} if
$f\left(N_i\right)\leq f\left(N_{i+1}\right)$, whenever $N_i$ and $N_{i+1}$ are vertically adjacent north steps. For $i=n$, this condition is imposed precisely when $P$ ends with a
north step, in which case it is interpreted using the above cyclic convention.
A pair $\pi=\left(P,f\right)$ is called a \emph{cyclic $\left(m,n\right)$-parking function} in the
strict case, and a \emph{weak cyclic $\left(m,n\right)$-parking function} in the weak case.
We denote the corresponding sets by $\cPF_{m,n}$ and
$\overline{\cPF}_{m,n}$. Every strictly cyclic labeling is weakly cyclic, i.e.
\(
 \cPF_{m,n}\subseteq\overline{\cPF}_{m,n}.
\)

For $k\geq1$, let $\cPF_{m,n}^k$ and
$\overline{\cPF}_{m,n}^k$ be the sets of ordered $k$-tuples
of elements of $\cPF_{m,n}$ and $\overline{\cPF}_{m,n}$,
respectively.
Write
$\pi^\bullet=(\pi^{(1)},\ldots,\pi^{(k)})$, with
$\pi^{(r)}=(P^{(r)},f^{(r)})$.
For such a tuple, we write
$P_{\pi^\bullet}=\left(P^{\left(1\right)},\ldots,P^{\left(k\right)}\right)$ for the underlying path tuple, and
$f_{\pi^\bullet}$ for the labeling on all north steps in all components.

Fix $\pi^\bullet\in\overline{\cPF}_{m,n}^k$.
A \emph{sheet} is a copy of the lattice strip $0\le y\le n$.
We use $k$ such copies, placing $\pi^{\left(r\right)}$ in sheet $r$ for $1\le r\le k$. If $c$ is a cell in the $r$-th sheet and $\left(x,y\right)$ is the
southeast corner of $c$, define its \emph{content} by
\[
    \operatorname{cont}\left(c\right)=k\left(my-nx\right)+r-1.
\]
The \emph{content} of a north step is the content of the cell immediately to
its left, and the content of an east step is the content of the cell immediately
below it.

Let $\pi$ be a cyclic parking function in $\cPF_{m,n}$ or $\overline{\cPF}_{m,n}$, and denote its $i$-th north step by $N_i$.
Consider a collection of indeterminates
$\left\{z_{i,j,r}\right\}$ indexed by
$1\le i\le n$, $j\in\mathbb Z_{\ge0}$, and $r\in\mathbb Z_{>0}$.
We define the \emph{skeleton} of $\pi$, denoted $z\left(\pi\right)$, as the monomial
\[
    z\left(\pi\right)=\prod_{i=1}^n
z_{i,\aseq\left(P_\pi\right)_i,f_\pi\left(N_i\right)}.
\]
For $\pi^\bullet=\left(\pi^{\left(1\right)},\dots,\pi^{\left(k\right)}\right)$ in either
$\cPF^k_{m,n}$ or $\overline{\cPF}^k_{m,n}$, we define its skeleton by
\[
    z\left(\pi^\bullet\right)=\prod_{\ell=1}^k z\left(\pi^{\left(\ell\right)}\right).
\]

The \emph{area} of a path $P$, denoted by $\area\left(P\right)$, is the number of cells below $P$ and above the
line $my=nx$, equivalently, the number of cells below $P$ with nonnegative contents. We have $\area\left(P\right)=\sum_{i}\aseq\left(P\right)_i$. For a tuple $\pi^\bullet$ in either
$\cPF^k_{m,n}$ or $\overline{\cPF}^k_{m,n}$, set
\[
    \area\left(P_{\pi^\bullet}\right)
    \coloneqq\sum_{r=1}^k \area\left(P^{\left(r\right)}\right),
    \qquad
    X_{\pi^\bullet}
    \coloneqq\prod_N x_{f_{\pi^\bullet}\left(N\right)},
\]
where the product is over all north steps in all components of
$P_{\pi^\bullet}$. We also write $\area\left(\pib\right)$ for $\area\left(P_{\pib}\right)$.

\begin{example}
For $m=3, n=2$, Figure~\ref{fig: 1} shows $\pi$ and $\pi'$ together with the line $my=nx$. We have $\pi\in \cPF_{m,n}$ (therefore also $\pi\in \overline{\cPF}_{m,n}$) since $\pi_1$ and $\pi_2$ are vertically adjacent north steps
whose labels are strictly increasing. Meanwhile, $\pi'\in \overline{\cPF}_{m,n}$ but
$\pi'\notin \cPF_{m,n}$, since $\pi'_2$ and
$\pi'_3$, the $\left(m,n\right)$-translate of $\pi'_1$, are vertically adjacent north steps whose labels are equal. We have
\begin{align*}
    \aseq\left(P_{\pi}\right)=\left(0,1\right), \quad
    \aseq\left(P_{\pi'}\right)=\left(2,0\right),\\
    z\left(\pi\right)=z_{1,0,1}z_{2,1,3}, \quad
    z\left(\pi'\right)=z_{1,2,4}z_{2,0,4}.
\end{align*}
    \begin{figure}[ht]
\centering

$\displaystyle \pi=$
\begin{tikzpicture}[scale=0.85, baseline=(current bounding box.center)]

    \draw[help lines, gray!30, thin] (-1,0) grid (3,2);
    \draw[black, thin] (0,0) -- (3,2);

    \draw[red, line width=2pt, line cap=round, line join=round]
        (0,0) -- (0,1) -- (0,2)
        -- (1,2) -- (2,2) -- (3,2);

    \node[fill=white, inner sep=1pt, text=red, font=\scriptsize]
        at (0,0.5) {$1$};
    \node[fill=white, inner sep=1pt, text=red, font=\scriptsize]
        at (0,1.5) {$3$};
\end{tikzpicture}
\qquad
$\displaystyle \pi'=$
\begin{tikzpicture}[scale=0.85, baseline=(current bounding box.center)]

    \draw[help lines, gray!30, thin] (-2,0) grid (3,2);
    \draw[black, thin] (0,0) -- (3,2);

    \draw[blue, line width=2pt, line cap=round, line join=round]
        (-2,0) -- (-2,1) -- (0,1)
        -- (1,1) -- (1,2) ;

    \node[fill=white, inner sep=1pt, text=blue, font=\scriptsize]
        at (-2,0.5) {$4$};
    \node[fill=white, inner sep=1pt, text=blue, font=\scriptsize]
        at (1,1.5) {$4$};
\end{tikzpicture}
\caption{Example of a cyclic parking function and a weak cyclic parking function for $m=3,n=2$.}
\label{fig: 1}
\end{figure}
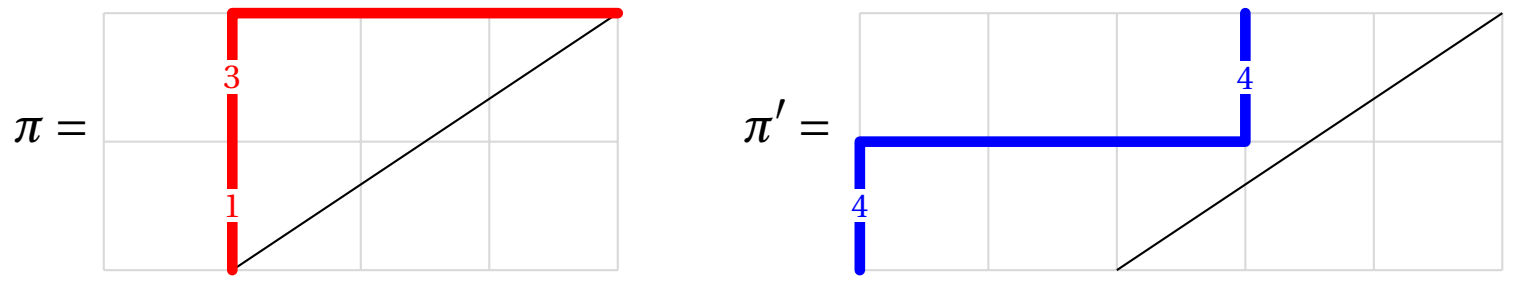
\end{example}

For $\pi^\bullet\in\overline{\cPF}_{m,n}^k$, the \emph{path diagonal inversion} statistic is
\[
    \pdinv\left(P_{\pi^\bullet}\right)
    \coloneqq
    \#\left\{\left(N,E\right): N \text{ is a north step},\ E \text{ is an east step},
    \ \operatorname{cont}\left(N\right)<\operatorname{cont}\left(E\right)\right\},
\]
where steps range over all components.

To define labeled diagonal inversions, first set
\[
    \mathcal A\left(P_{\pi^\bullet}\right)
    \coloneqq 
    \left\{\left(N,N'\right): \operatorname{cont}\left(N\right)<\operatorname{cont}\left(N'\right)
    <\operatorname{cont}\left(N\right)+km\right\}
\]
for $\pi^\bullet\in\overline{\cPF}_{m,n}^k$.
Its elements are called \emph{attacking pairs}.
Here $N$ and $N'$ range over all north steps in all components
of $P_{\pi^\bullet}$. For $\pi^\bullet\in\cPF^k_{m,n}$, define
\[
    \ldinv\left(\pi^\bullet\right)
\coloneqq 
\#\left\{\left(N,N'\right)\in \mathcal A\left(P_{\pi^\bullet}\right):
f_{\pi^\bullet}\left(N\right)\geq f_{\pi^\bullet}\left(N'\right)\right\}.
\]
For $\pi^\bullet\in\overline{\cPF}^k_{m,n}$, define its weak analogue by
\[
    \overline{\ldinv}\left(\pi^\bullet\right)
\coloneqq 
\#\left\{\left(N,N'\right)\in \mathcal A\left(P_{\pi^\bullet}\right):
f_{\pi^\bullet}\left(N\right)> f_{\pi^\bullet}\left(N'\right)\right\}.
\]

Finally, define
\[
    C\left(m,n,k\right)\coloneqq\frac{\left(mk-1\right)\left(nk-1\right)+k-1}{2}.
\]
For $\pi^\bullet\in\cPF^k_{m,n}$, set
\[
    \stat\left(\pi^\bullet\right)
    \coloneqq C\left(m,n,k\right)-\pdinv\left(P_{\pi^\bullet}\right)-\ldinv\left(\pi^\bullet\right),
\]
while for $\pi^\bullet\in\overline{\cPF}^k_{m,n}$, set
\[
    \overline{\stat}\left(\pi^\bullet\right)
    \coloneqq C\left(m,n,k\right)-\pdinv\left(P_{\pi^\bullet}\right)-\overline{\ldinv}\left(\pi^\bullet\right).
\]

A cyclic parking function is called \emph{standard} if its labels are
$1,2,\ldots,n$, each appearing exactly once. More generally, a labeling of a
$k$-tuple is called \emph{standard} if $f_{\pi^\bullet}$ is a bijection from the
set of all north steps of $P_{\pi^\bullet}$ to $\left\{1,2,\ldots,kn\right\}$. A weak cyclic tuple with standard labels is strictly cyclic, and
$\ldinv=\overline{\ldinv}$ on such tuples. We write
$\cPF^{k,\mathrm{std}}_{m,n}$ for the subset of tuples with standard labels.

We use two standardization maps, one for each convention. The strict
standardization $\operatorname{std}$, defined on $\cPF^k_{m,n}$, keeps the
underlying path tuple fixed and relabels all north steps by $1,2,\ldots,kn$ in
increasing order of their original labels, breaking ties by decreasing content.
The weak standardization $\overline{\operatorname{std}}$, defined on
$\overline{\cPF}^k_{m,n}$, is defined in the same way, except that ties are
broken by increasing content. Any remaining ties are broken by a fixed
auxiliary order. With these conventions,
\[
    \ldinv\left(\operatorname{std}\left(\pi^\bullet\right)\right)=\ldinv\left(\pi^\bullet\right),
    \qquad
    \overline{\ldinv}\left(\overline{\operatorname{std}}\left(\pi^\bullet\right)\right)
    =\overline{\ldinv}\left(\pi^\bullet\right).
\]
Consequently, $\stat$ and $\overline{\stat}$ are preserved under their
respective standardizations.

For labeled north-step types at the same height, we use the total order
\begin{equation}\label{eq:labeled-step-order}
 \left(N,a\right)<\left(N',a'\right)
 \quad\Longleftrightarrow\quad
 x\left(N\right)>x\left(N'\right)\ \text{or}\ \left(x\left(N\right)=x\left(N'\right)\text{ and }a>a'\right),
\end{equation}
where $x\left(N\right)$ is the horizontal coordinate of $N$.
Thus the order reads positions from right to left, breaking ties by
decreasing labels. 

Let $\left(N, a\right)$ and $\left(N', a'\right)$ be two labeled north steps located at consecutive heights. Let $\left(x, y\right)$ be the endpoint of $N$ and $\left(x', y\right)$ be the starting point of $N'$. We say that $\left(N,a\right)$ can be \emph{extended}
(respectively, \emph{weakly extended}) with $\left(N',a'\right)$
if $x<x'$, or if $x=x'$ and $a<a'$
(respectively, $a\leq a'$).
By definition, \(\pi\in\cPF_{m,n}\) (respectively, \(\pi\in\overline{\cPF}_{m,n}\)) if and only if \(\pi_\ell\) can be extended (respectively, weakly extended) with \(\pi_{\ell+1}\) for every \(1\le\ell\le n\).

For $1\le \ell\le n$ and $1\le i<j\le k$, we say that
$\pi^\bullet\in\cPF_{m,n}^k$ has a \emph{crossing} at
$(\ell,i,j)$ if $\pi^{(i)}_\ell$ can be extended with
$\pi^{(j)}_{\ell+1}$ and $\pi^{(j)}_\ell$ can be extended
with $\pi^{(i)}_{\ell+1}$. A crossing is called
\emph{positive} if
$\pi^{(i)}_\ell>\pi^{(j)}_\ell$ in the order defined above.
For $\pi^\bullet\in\overline{\cPF}_{m,n}^k$, we define
\emph{weak crossings} and \emph{weak positive crossings}
analogously, replacing both occurrences of
``can be extended'' by ``can be weakly extended''.
We write $\pos(\pi^\bullet)$ and
$\overline{\pos}(\pi^\bullet)$ for the numbers of positive
crossings and weak positive crossings, respectively.
Heights are interpreted cyclically: when $\ell=n$, the
north steps at height $1$ are translated by $(m,n)$.
Equal labeled north-step types do not form a positive
crossing or a weak positive crossing.

The following lemma relates these crossing statistics
to $\stat$ and $\overline{\stat}$.

\begin{lem}\label{lem: pos crossing}
The following statements hold.
\begin{enumerate}
    \item For $\pi^\bullet\in\cPF_{m,n}^k$, the difference
    \[
        \stat(\pi^\bullet)-\pos(\pi^\bullet)
    \]
    depends only on the skeleton $z(\pi^\bullet)$.

    \item For $\pi^\bullet\in\overline{\cPF}_{m,n}^k$,
    the difference
    \[
        \overline{\stat}(\pi^\bullet)
        -\overline{\pos}(\pi^\bullet)
    \]
    depends only on the skeleton $z(\pi^\bullet)$.
\end{enumerate}
\end{lem}

\begin{proof}
Part~(1) is precisely \cite[Lemma~3.2]{KO26}.
The proof of part~(2) is entirely analogous and is omitted.
\end{proof}

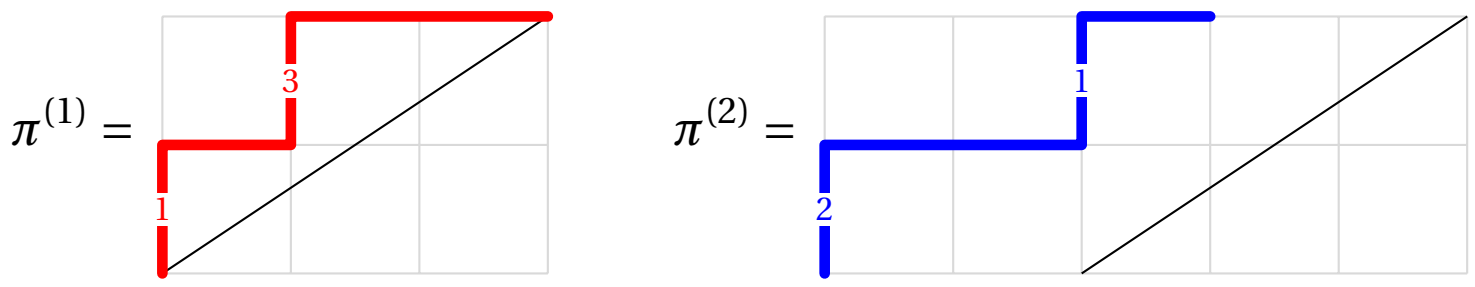
\begin{figure}[ht]
\centering

$\displaystyle \pi^{\left(1\right)}=$
\begin{tikzpicture}[scale=0.85, baseline=(current bounding box.center)]
    
    \draw[help lines, gray!30, thin] (0,0) grid (3,2);
    \draw[black, thin] (0,0) -- (3,2);

    \draw[red, line width=2pt, line cap=round, line join=round]
        (0,0) -- (0,1) -- (1,1)
        -- (1,2) -- (2,2) -- (3,2);

    \node[fill=white, inner sep=1pt, text=red, font=\scriptsize]
        at (0,0.5) {$1$};
    \node[fill=white, inner sep=1pt, text=red, font=\scriptsize]
        at (1,1.5) {$3$};
\end{tikzpicture}
\qquad
$\displaystyle \pi^{\left(2\right)}=$
\begin{tikzpicture}[scale=0.85, baseline=(current bounding box.center)]

    \draw[help lines, gray!30, thin] (-2,0) grid (3,2);
    \draw[black, thin] (0,0) -- (3,2);

    \draw[blue, line width=2pt, line cap=round, line join=round]
        (-2,0) -- (-2,1) -- (0,1)
        -- (0,2) -- (1,2);

    \node[fill=white, inner sep=1pt, text=blue, font=\scriptsize]
        at (-2,0.5) {$2$};
    \node[fill=white, inner sep=1pt, text=blue, font=\scriptsize]
        at (0,1.5) {$1$};
\end{tikzpicture}

\caption{Weak cyclic parking functions in $\overline{\cPF}_{3,2}$.}
\label{fig: crossing}
\end{figure}

\begin{example}
    Let $\tau=\left(\pi^{\left(1\right)},\pi^{\left(2\right)}\right)\in \overline{\cPF}_{3,2}^{2}$ and $\rho=\left(\pi^{\left(2\right)},\pi^{\left(1\right)}\right)\in \overline{\cPF}_{3,2}^{2}$ for $\pi^{\left(1\right)}$ and $\pi^{\left(2\right)}$ in Figure~\ref{fig: crossing}. Both $\tau$ and $\rho$ have a weak crossing at $\left(1,1,2\right)$, and this weak crossing is positive only for $\rho$. Therefore $\overline{\pos}\left(\tau\right)=0$ and $\overline{\pos}\left(\rho\right)=1$. We have
    $\pdinv\left(P_{\tau}\right)=5$,  $\pdinv\left(P_{\rho}\right)=3$, $\overline{\ldinv}\left(\tau\right)=1$, and $\overline{\ldinv}\left(\rho\right)=2$. We can check that
    \begin{equation*}
        \overline{\stat}\left(\tau\right)-\overline{\pos}\left(\tau\right)=\overline{\stat}\left(\rho\right)-\overline{\pos}\left(\rho\right).
    \end{equation*}
\end{example}

\section{\texorpdfstring{$\h$}{h}-operators}\label{Sec: h operators}
In this section, we give a combinatorial description of
$\omega_X\left(e_k[-X^{m,n}Y]\cdot 1\right)$; see Proposition~\ref{prop: omega Omega}.

Fix coprime positive integers $m,n$.
We use two strict partial orders:
$\prec^{\mathrm s}=\prec^{\mathrm s}_{m,n}$ on $\cPF_{m,n}$
and $\prec^{\mathrm w}=\prec^{\mathrm w}_{m,n}$
on $\overline{\cPF}_{m,n}$.
For $\pi=(P,f)$ and $\pi'=(P',f')$ in $\cPF_{m,n}$,
we write $\pi\prec^{\mathrm s}\pi'$ if, for every
$1\le\ell\le n$, the labeled north step $\pi_\ell$
cannot be extended with $\pi'_{\ell+1}$.
In this case, $P_{\pi'}$ lies weakly to the left of $P_\pi$.
Similarly, for $\pi,\pi'\in\overline{\cPF}_{m,n}$,
we write $\pi\prec^{\mathrm w}\pi'$ if, for every
$1\le\ell\le n$, $\pi_\ell$ cannot be weakly extended
with $\pi'_{\ell+1}$.

The pair $\left(\pi,\pi'\right)\in\cPF_{m,n}^{2}$ has no crossing if and only if
$\pi\prec^{\mathrm{s}}\pi'$ or $\pi'\prec^{\mathrm{s}}\pi$.
An analogous statement holds for pairs in $\overline{\cPF}_{m,n}^{2}$.

\subsection{\texorpdfstring{$\mathfrak h$}{h}-operators}
\label{subsec:h-operators}

For $a\geq 1$, let $\CH_a$ be the set of length-$a$ chains
$\left(\pi^{\left(1\right)},\ldots,\pi^{\left(a\right)}\right)$ such that
$\pi^{\left(i\right)}\in\cPF_{m,n}$ for $1\leq i\leq a$ and
\[
    \pi^{\left(i\right)}\prec^{\mathrm{s}}\pi^{\left(i+1\right)}
    \qquad \left(1\leq i<a\right).
\]
Similarly, let $\overline{\CH}_a$ be the set of length-$a$
weak chains $\Pi=\left(\pi^{\left(1\right)},\ldots,\pi^{\left(a\right)}\right)$ such that
$\pi^{\left(i\right)}\in\overline{\cPF}_{m,n}$ for $1\leq i\leq a$ and
\[
    \pi^{\left(i\right)}\prec^{\mathrm{w}}\pi^{\left(i+1\right)}
    \qquad \left(1\leq i<a\right).
\]
We also set $\CH_0=\overline{\CH}_0=\left\{\emptyset\right\}$.

We recall the $\mathfrak h$-operators introduced in
\cite[Section~4]{KO26}, which describe symmetric functions
of the form $f\left[-MX^{m,n}\right]\cdot 1$ in terms of generating
functions for chains in $\cPF_{m,n}$.
We also introduce their weak analogues, the
$\mathfrak h^{\mathrm w}$-operators, to describe
$\omega_X\!\left(f\left[-MX^{m,n}\right]\cdot 1\right)$.

Let
\[
    \mathcal D \subset
    \mathbb C\left(q,t\right)\left[\!\left[z_{i,j,r}:
        1\leq i\leq n,\ j\in\mathbb Z_{\geq0},\
        r\in\mathbb Z_{>0}
    \right]\!\right]
\]
be the completed subalgebra generated by the skeleton
monomials $z\left(\pi\right)$ for $\pi\in\cPF_{m,n}$.
Following \cite[Definition~4.2]{KO26}, we define
$\mathfrak h_a$ on a skeleton monomial $z\left(\tau^\bullet\right)$ by
\[
    \mathfrak h_a\cdot z\left(\tau^\bullet\right)
    \coloneqq
    \sum_{\pi^\bullet\in\CH_a}
    q^{\stat\left(\tau^\bullet,\pi^\bullet\right)-\stat\left(\tau^\bullet\right)}
    z\left(\tau^\bullet,\pi^\bullet\right),
\]
where $\left(\tau^\bullet,\pi^\bullet\right)$ denotes concatenation.
The action is extended linearly and continuously to
$\mathcal D$.
The right-hand side depends only on $z\left(\tau^\bullet\right)$,
as verified in \cite[Definition~4.2]{KO26}.
We use the conventions $z\left(\emptyset\right)=1$ and
$\stat\left(\emptyset\right)=0$; in particular, $\mathfrak h_0$ is
the identity operator.

For a weak composition $\alpha=\left(\alpha_1,\ldots,\alpha_\ell\right)$,
write
\[
    \mathfrak h_\alpha
    \coloneqq
    \mathfrak h_{\alpha_\ell}\cdots\mathfrak h_{\alpha_1},
    \qquad
    \CH_\alpha
    \coloneqq
    \CH_{\alpha_1}\times\cdots\times\CH_{\alpha_\ell}.
\]
We regard an element of $\CH_\alpha$ as a single tuple
by concatenation.
Iterating the defining action gives
\[
    \mathfrak h_\alpha\cdot 1
    =
    \sum_{\pi^\bullet\in\CH_\alpha}
    q^{\stat\left(\pi^\bullet\right)}z\left(\pi^\bullet\right).
\]
The operators commute:
\[
    \mathfrak h_a\mathfrak h_b
    =
    \mathfrak h_b\mathfrak h_a
    \qquad \left(a,b\geq0\right),
\]
as proved in \cite[Equation~(4.1)]{KO26}.

We now specialize $z_{i,j,r}\mapsto t^j x_r$.
Under this specialization, $z\left(\pi^\bullet\right)$ becomes
$t^{\area\left(\pi^\bullet\right)}X_{\pi^\bullet}$.
By abuse of notation, we continue to write
$\mathfrak h_\alpha\cdot 1$ for the specialized series:
\[
    \mathfrak h_\alpha\cdot 1
    =
    \sum_{\pi^\bullet\in\CH_\alpha}
    t^{\area\left(\pi^\bullet\right)}
    q^{\stat\left(\pi^\bullet\right)}X_{\pi^\bullet}.
\]
Throughout, this specialization is applied after the
operators have acted on $1$.
The following result relates these series to the elliptic
Hall algebra action and, in particular, implies that they
are symmetric in the $X$-variables.

\begin{thm}[{\cite[Corollary~5.14]{KO26}}]
\label{thm: ko 26 main}
For any symmetric function $f$, we have
\[
    f\left[-MX^{m,n}\right]\cdot 1
    =
    \sum_\lambda c_\lambda
    \left(\mathfrak h_\lambda\cdot 1\right),
\]
where $c_\lambda$ is the coefficient of $h_\lambda$ in
the expansion of $f\left[\left(1-t\right)X\right]$.
\end{thm}

To obtain an analogous description of
$\omega_X\!\left(f\left[-MX^{m,n}\right]\cdot 1\right)$, we introduce
weak $\mathfrak h$-operators.
Let $\overline{\mathcal D}$ be the completed subalgebra
of the same power-series ring generated by $z\left(\pi\right)$ for
$\pi\in\overline{\cPF}_{m,n}$.
We define
\[
    \mathfrak h_a^{\mathrm w}\cdot z\left(\tau^\bullet\right)
    \coloneqq
    \sum_{\pi^\bullet\in\overline{\CH}_a}
    q^{\overline{\stat}\left(\tau^\bullet,\pi^\bullet\right)
       -\overline{\stat}\left(\tau^\bullet\right)}
    z\left(\tau^\bullet,\pi^\bullet\right),
\]
and extend the action linearly and continuously to
$\overline{\mathcal D}$.
As in the strict case, well-definedness follows from the
skeleton invariance of
$\overline{\stat}-\overline{\pos}$
in Lemma~\ref{lem: pos crossing}.
We set $\overline{\stat}\left(\emptyset\right)=0$, so that
$\mathfrak h_0^{\mathrm w}$ is the identity operator.

For a weak composition $\alpha=\left(\alpha_1,\ldots,\alpha_\ell\right)$,
write
\[
    \mathfrak h_\alpha^{\mathrm w}
    \coloneqq
    \mathfrak h_{\alpha_\ell}^{\mathrm w}\cdots
    \mathfrak h_{\alpha_1}^{\mathrm w},
    \qquad
    \overline{\CH}_\alpha
    \coloneqq
    \overline{\CH}_{\alpha_1}\times\cdots
    \times\overline{\CH}_{\alpha_\ell},
\]
again identifying tuples by concatenation.
The defining action gives
\[
    \mathfrak h_\alpha^{\mathrm w}\cdot 1
    =
    \sum_{\pi^\bullet\in\overline{\CH}_\alpha}
    q^{\overline{\stat}\left(\pi^\bullet\right)}z\left(\pi^\bullet\right).
\]
The weak operators satisfy the corresponding commutativity
relation:
\[
    \mathfrak h_a^{\mathrm w}\mathfrak h_b^{\mathrm w}
    =
    \mathfrak h_b^{\mathrm w}\mathfrak h_a^{\mathrm w}
    \qquad \left(a,b\geq0\right).
\]
The proof is entirely analogous to the proof of
commutativity for the strict operators in \cite[Equation~(4.1)]{KO26},
so we omit it.

Using the same specialization and notation convention,
we write
\[
    \mathfrak h_\alpha^{\mathrm w}\cdot 1
    =
    \sum_{\pi^\bullet\in\overline{\CH}_\alpha}
    t^{\area\left(\pi^\bullet\right)}
    q^{\overline{\stat}\left(\pi^\bullet\right)}X_{\pi^\bullet}.
\]
The following proposition relates these series to their
strict counterparts and implies their symmetry in $X$.

\begin{proposition}\label{prop:omega-h-hw}
For a partition $\lambda$, we have
\[
    \omega\left(\mathfrak h_\lambda\cdot 1\right)
    =
    \mathfrak h^{\mathrm{w}}_\lambda\cdot 1.
\]
\end{proposition}

\begin{proof}
Write $\lambda=\left(\lambda_1,\ldots,\lambda_\ell\right)$ and set $k=\left|\lambda\right|$.
We regard $\CH_\lambda$ and $\overline{\CH}_\lambda$ as subsets of
$\cPF_{m,n}^{k}$ and $\overline{\cPF}_{m,n}^{k}$, respectively, by
concatenating the chains. Let
\[
    \CH_\lambda^{\mathrm{std}}
    \coloneqq
    \CH_\lambda\cap \cPF_{m,n}^{k,\mathrm{std}}.
\]

Fix \(\pi^\bullet\in\CH_\lambda^{\mathrm{std}}\). For each \(i\in[kn]\), let \(N(i)\coloneqq f_{\pi^\bullet}^{-1}(i)\) be the north step carrying label \(i\). Define
\[
    D\left(\pi^\bullet\right)
    \coloneqq
    \left\{
        i\in\left[kn-1\right]:
        \cont\left(N\left(i\right)\right)<\cont\left(N\left(i+1\right)\right)
    \right\}.
\]
We write $D\left(\pi^\bullet\right)^c$ for its complement in $\left[kn-1\right]$.

Now define
\[
    \CH_\lambda\left(\pi^\bullet\right)
    \coloneqq
    \left\{
        \rho^\bullet\in\CH_\lambda:
        \operatorname{std}\left(\rho^\bullet\right)=\pi^\bullet
    \right\}
\qquad \text{and} \qquad    \overline{\CH}_\lambda\left(\pi^\bullet\right)
    \coloneqq
    \left\{
        \rho^\bullet\in\overline{\CH}_\lambda:
        \overline{\operatorname{std}}\left(\rho^\bullet\right)=\pi^\bullet
    \right\}.
\]
For $\rho^\bullet$ in either
$\CH_\lambda\left(\pi^\bullet\right)$ or
$\overline{\CH}_\lambda\left(\pi^\bullet\right)$, set
\[
 f_i\coloneqq f_{\rho^\bullet}\left(N\left(i\right)\right)
 \qquad \left(1\le i\le kn\right).
\]
For $\rho^\bullet\in\CH_\lambda(\pi^\bullet)$,
this label sequence satisfies
\begin{equation}\label{eq:strict-labeling-condition}
    f_1\leq\cdots\leq f_{kn},
    \qquad
    f_i<f_{i+1}\quad\text{for }i\in D\left(\pi^\bullet\right).
\end{equation}
Similarly, for
$\rho^\bullet\in\overline{\CH}_\lambda\left(\pi^\bullet\right)$,
the corresponding label sequence satisfies
\begin{equation}\label{eq:weak-labeling-condition}
    f_1\leq\cdots\leq f_{kn},
    \qquad
    f_i<f_{i+1}\quad\text{for }i\in D\left(\pi^\bullet\right)^c.
\end{equation}

Conversely, any labeling of the underlying path tuple $P_{\pi^\bullet}$
satisfying \eqref{eq:strict-labeling-condition}
(respectively, \eqref{eq:weak-labeling-condition})
gives an element of $\CH_\lambda\left(\pi^\bullet\right)$
(respectively, $\overline{\CH}_\lambda\left(\pi^\bullet\right)$). Therefore we conclude
\[
\begin{aligned}
    \sum_{\rho^\bullet\in\CH_\lambda\left(\pi^\bullet\right)}
    q^{\stat\left(\rho^\bullet\right)}
    t^{\area\left(P_{\rho^\bullet}\right)}
    X_{\rho^\bullet}
    &=
    q^{\stat\left(\pi^\bullet\right)}
    t^{\area\left(P_{\pi^\bullet}\right)}
    F_{D\left(\pi^\bullet\right)}\left[X\right],\\
    \sum_{\rho^\bullet\in
        \overline{\CH}_\lambda\left(\pi^\bullet\right)}
    q^{\overline{\stat}\left(\rho^\bullet\right)}
    t^{\area\left(P_{\rho^\bullet}\right)}
    X_{\rho^\bullet}
    &=
    q^{\stat\left(\pi^\bullet\right)}
    t^{\area\left(P_{\pi^\bullet}\right)}
    F_{D\left(\pi^\bullet\right)^c}\left[X\right].
\end{aligned}
\]
Here, $F_E\left[X\right]$ is Gessel's \emph{fundamental quasisymmetric function} \cite{Ges84} for $E\subseteq\left[kn-1\right]$, defined by
\[
    F_E\left[X\right]
    \coloneqq
    \sum_{\substack{
        1\leq i_1\leq\cdots\leq i_{kn}\\
        j\in E\Rightarrow i_j<i_{j+1}
    }}
    x_{i_1}\cdots x_{i_{kn}}.
\]
We suppress the degree in the notation for fundamental
quasisymmetric functions. In degree $d$, complements of
indexing subsets are taken in $\left[d-1\right]$.
Then the involution $\omega$ on $\Lambda$ can be extended to the ring of quasisymmetric functions by setting
\begin{equation}\label{eq: involution on F}
    \omega \left(F_{E}\right)=F_{E^c},
\end{equation} 
where $E^c=[d-1]\setminus E$ when $F_E$ has degree $d$
(see, for example, \cite[Section~3.6]{LMvW13}).

Therefore,
\[
    \omega\left(
    \sum_{\rho^\bullet\in\CH_\lambda\left(\pi^\bullet\right)}
    q^{\stat\left(\rho^\bullet\right)}
    t^{\area\left(P_{\rho^\bullet}\right)}
    X_{\rho^\bullet}
    \right)
    =
    \sum_{\rho^\bullet\in\overline{\CH}_\lambda\left(\pi^\bullet\right)}
    q^{\overline{\stat}\left(\rho^\bullet\right)}
    t^{\area\left(P_{\rho^\bullet}\right)}
    X_{\rho^\bullet}.
\]
Summing over all
$\pi^\bullet\in\CH_\lambda^{\mathrm{std}}$ proves the proposition.
\end{proof}

\subsection{\texorpdfstring{$\mathfrak e$}{e}-operators}
\label{subsec:e-operators}

Since the operators $\mathfrak h^{\mathrm{w}}_a$ commute, the following
Jacobi--Trudi determinant is well-defined. We take \(\mathfrak h^{\mathrm w}_0\) to be the identity operator and set \(\mathfrak h^{\mathrm w}_r=0\) for \(r<0\). For
$a\geq 1$, define
\begin{equation*}
    \mathfrak e^{\mathrm{w}}_a
    =
    \det\left(\mathfrak h^{\mathrm{w}}_{1-i+j}\right)_{1\leq i,j\leq a}.
\end{equation*}
Equivalently,
\begin{equation}\label{eq: e to h}
    \mathfrak e^{\mathrm{w}}_a
    =
    \det
    \begin{pmatrix}
        \mathfrak h^{\mathrm{w}}_1
            & \mathfrak h^{\mathrm{w}}_2
            & \cdots
            & \mathfrak h^{\mathrm{w}}_a \\
        1
            & \mathfrak h^{\mathrm{w}}_1
            & \cdots
            & \mathfrak h^{\mathrm{w}}_{a-1} \\
        0
            & 1
            & \cdots
            & \mathfrak h^{\mathrm{w}}_{a-2} \\
        \vdots
            & \vdots
            & \ddots
            & \vdots \\
        0
            & 0
            & \cdots
            & \mathfrak h^{\mathrm{w}}_1
    \end{pmatrix}.
\end{equation}
Let $\overline{\CE}_a$ be the set of tuples
$\left(\pi^{\left(1\right)},\ldots,\pi^{\left(a\right)}\right)$ such that
$\pi^{\left(i\right)}\in\overline{\cPF}_{m,n}$ and
\[
    \pi^{\left(i\right)}\nprec^{\mathrm w}\pi^{\left(i+1\right)}
    \qquad\text{for }1\leq i\leq a-1.
\]
For a composition $\alpha=\left(\alpha_1,\ldots,\alpha_\ell\right)$, define
\[
    \mathfrak e^{\mathrm w}_\alpha
    \coloneqq
    \mathfrak e^{\mathrm w}_{\alpha_\ell}\cdots
    \mathfrak e^{\mathrm w}_{\alpha_1},
    \qquad
    \overline{\CE}_\alpha
    \coloneqq
    \overline{\CE}_{\alpha_1}\times\cdots\times
    \overline{\CE}_{\alpha_\ell}.
\]
An inclusion--exclusion argument applied to \eqref{eq: e to h} gives
\[
    \mathfrak e^{\mathrm w}_\alpha\cdot 1
    =
    \sum_{\pi^\bullet\in\overline{\CE}_\alpha}
    q^{\overline{\stat}\left(\pi^\bullet\right)}z\left(\pi^\bullet\right).
\]

As before, we use the same notation after specializing
$z_{i,j,r}\mapsto t^j x_r$.
The following proposition gives the operator expansion that
will be compared with the affine Springer fiber formula
in the next section.

\begin{proposition}\label{prop: omega Omega}
We have
\begin{equation}\label{eq:omega-Omega-ew}
    \omega_X\!\left(e_k[-X^{m,n}Y]\cdot1\right)
    =
    \sum_{\mu\vdash k}
    \left(\mathfrak e^{\mathrm{w}}_\mu\cdot 1\right)
    m_\mu\!\left[\frac{Y}{1-q}\right].
\end{equation}
\end{proposition}

\begin{proof}
For partitions \(\lambda,\mu\vdash k\), let $c_{\lambda,\mu}$ (respectively, $c'_{\lambda,\mu}$) be the coefficient of $h_{\mu}$ (respectively, $e_{\mu}$) in $e_{\lambda}\left[\frac{X}{1-q}\right]$. Theorem~\ref{thm: ko 26 main} and Proposition~\ref{prop:omega-h-hw} give
\begin{equation*}
   \omega_X\left( e_{\lambda}\left[-X^{m,n}\right]\cdot 1\right)=\sum_{\mu \vdash k} c_{\lambda,\mu} \left(\h^{\mathrm{w}}_{\mu}\cdot 1\right).
\end{equation*}

Since the operators $\mathfrak h_a^{\mathrm w}$ commute,
$\mathfrak e_\mu^{\mathrm w}$ is obtained from the usual
expansion of $e_\mu$ in the $h$-basis by replacing each
$h_\nu$ with $\mathfrak h_\nu^{\mathrm w}$.  Hence
\[
    \omega_X\left(e_\lambda\left[-X^{m,n}\right]\cdot 1\right)
    =
    \sum_{\mu\vdash k}
    c'_{\lambda,\mu}\,
    \mathfrak e^{\mathrm w}_\mu\cdot 1.
\]
By the dual Cauchy identity, we have
\[
    \sum_{\nu\vdash k}
    e_\nu\!\left[\frac{X}{1-q}\right]m_\nu\left[Y\right]
    =
    e_k\!\left[\frac{XY}{1-q}\right]
    =
    \sum_{\nu\vdash k}
    e_\nu\left[X\right]m_\nu\!\left[\frac{Y}{1-q}\right].
\]
Comparing the coefficients of $e_\mu\left[X\right]m_\lambda\left[Y\right]$ gives
\begin{equation}\label{eq: e m}
    c'_{\lambda,\mu}
    =
    \left[e_\mu\left[X\right]\right]\,e_\lambda\!\left[\frac{X}{1-q}\right]
    =
    \left[m_\lambda\left[Y\right]\right]\,m_\mu\!\left[\frac{Y}{1-q}\right].
\end{equation}
Therefore we have
\begin{align*}
\omega_X\left(e_k[-X^{m,n}Y]\cdot1\right)
&=
\sum_{\lambda\vdash k}
\omega_X\left(e_\lambda\left[-X^{m,n}\right]\cdot1\right)m_\lambda\left[Y\right] \\
&=
\sum_{\lambda,\mu\vdash k}
c'_{\lambda,\mu}
\left(\mathfrak e^{\mathrm w}_\mu\cdot1\right)m_\lambda\left[Y\right] \\
&=
\sum_{\mu\vdash k}
\left(\mathfrak e^{\mathrm w}_\mu\cdot1\right)
m_\mu\!\left[\frac{Y}{1-q}\right].
\end{align*}

\end{proof}

\section{Proof of Theorem~\ref{thm: main}}\label{Sec: asf combinatorics}

\subsection{Cyclic chromatic symmetric functions}
A Hessenberg function of size $k$ is a weakly increasing sequence
$\left(a_1,\dots,a_k\right)$ of integers such that
\[
    i\le a_i\le k
    \qquad\text{for all }1\le i\le k.
\]
A Hessenberg function of size \(k\) is often identified with a Dyck path of size \(k\), namely a lattice path from \(\left(0,0\right)\) to \(\left(k,k\right)\) staying weakly above the diagonal \(y=x\).

 We recall the cyclic chromatic symmetric functions introduced by Boixeda Alvarez--Hameister--Ma~\cite{BAHM26}. Let $\mm^{\bullet}=\left(\mm^{\left(1\right)},\mm^{\left(2\right)},\dots,\mm^{\left(n\right)}\right)$ be an $n$-tuple of Hessenberg functions of size $k$. We denote the set of such tuples by \(\Hessen\left(k,n\right)\).

To each \(\mm^\bullet\in\Hessen(k,n)\), we associate a directed graph \(G(\mm^\bullet)\) with vertex set
\begin{equation*}
    V\left(G\left(\mm^{\bullet}\right)\right)=\left\{\left(i,\ell\right): 1 \leq i \leq k, 1 \leq \ell \leq n\right\}.
\end{equation*}
We regard the second coordinate cyclically, identifying \(\left(i,n+1\right)\) with \(\left(i,1\right)\). Directed edges are given as follows: there exists an edge from $\left(i,\ell\right)$ to $\left(j,\ell+1\right)$ when $i\leq \mm^{\left(\ell\right)}_{j}$.

A \emph{loop} \(L\) of \(\mm^\bullet\) is a subset of the vertex set of \(G\left(\mm^\bullet\right)\) of the form
\begin{equation*}
    \left\{\left(a_1,1\right),\left(a_2,2\right),\dots,\left(a_n,n\right)\right\}
\end{equation*}
such that there exists a directed edge from $\left(a_{\ell},\ell\right)$ to $\left(a_{\ell+1},\ell+1\right)$ for every \(1\le\ell\le n\), with indices interpreted modulo \(n\). We write $L\left(\ell\right)=a_\ell$ for $1\leq\ell\leq n$. A \emph{loop decomposition} of \(\mm^\bullet\) is an ordered
\(k\)-tuple \(\left(L_1,\ldots,L_k\right)\) of pairwise disjoint loops of
\(\mm^\bullet\) whose union is \(V\left(G\left(\mm^\bullet\right)\right)\).
The set of loop decompositions of \(G\left(\mm^\bullet\right)\) is denoted by \(\LD\left(\mm^\bullet\right)\).

 For two loops $L_1,L_2$ of $\mm^\bullet$, we write 
\begin{equation*}
    L_1\prec L_2
\end{equation*}
if there is no edge from a vertex in $L_1$ to a vertex in $L_2$. Note that \(L_1\prec L_2\) implies \(L_1(\ell)>L_2(\ell)\) for every \(1\le\ell\le n\).

For a composition $\beta\models k$ and a loop decomposition $L=\left(L_1,L_2,\dots,L_k\right)\in \LD\left(\mm^{\bullet}\right)$,
\begin{itemize}
\item we say that $L$ is \emph{$\beta$-suitable} if
\begin{equation*}
    L_{r}\nprec L_{r+1}\nprec \dots\nprec L_{r+\beta_i-1}
\end{equation*}
for each \(1\le i\le\ell(\beta)\), where \(r=1+\sum_{j=1}^{i-1}\beta_j\).
    \item we say that $L$ is \emph{$\beta$-colorable} if
\begin{equation*}
    L_{r}\prec L_{r+1}\prec \dots\prec L_{r+\beta_i-1}
\end{equation*}
for each \(1\le i\le\ell(\beta)\), where \(r=1+\sum_{j=1}^{i-1}\beta_j\).
\end{itemize}
We denote the set of $\beta$-suitable
(respectively, $\beta$-colorable) loop decompositions
by $\LD\left(\mm^\bullet;\beta\right)$
(respectively, $\overline{\LD}\left(\mm^\bullet;\beta\right)$).

We say that a loop decomposition $L=\left(L_1,L_{2},\dots,L_k\right)$ has an \emph{ascent} at $\left(\ell,i,j\right)$ for some $\ell$ and $i<j$ if 
\begin{equation*}
    L_i\left(\ell\right)<L_j\left(\ell\right)\leq\mm^{\left(\ell\right)}_{L_i\left(\ell+1\right)},
\end{equation*}
where $\ell,\ell+1$ are taken modulo $n$. We denote the number of ascents of $L$ by $\asc(L)$.

Finally, we define the \emph{cyclic chromatic symmetric function}
$\Xi_{\mm^\bullet}\left[X;q\right]$
and the \emph{dual cyclic chromatic symmetric function}
$\widetilde{\Xi}_{\mm^\bullet}\left[X;q\right]$ by \footnote{Our ascent statistic is complementary to the
descent statistic in \cite{BAHM26}: their sum is
$D_{\mmb}=\sum_{\ell=1}^n\sum_{j=1}^k
(\mm_j^{(\ell)}-j)$.
Thus the two definitions are related by replacing $q$
with $q^{-1}$ and multiplying by $q^{D_{\mmb}}$; our definition of $\Xi_{\mmb}$ therefore agrees with
  theirs by~\cite[Proposition~3.10]{BAHM26}.}
\begin{align*}
    \Xi_{\mm^{\bullet}}\left[X;q\right]&=\sum_{\lambda\vdash k}\sum_{L\in \LD\left(\mm^{\bullet};\lambda\right)} q^{\asc\left(L\right)}m_{\lambda}\left[X\right],\\
    \widetilde{\Xi}_{\mm^{\bullet}}\left[X;q\right]&=\sum_{\lambda\vdash k}\sum_{L\in \overline{\LD}\left(\mm^{\bullet};\lambda\right)} q^{\asc\left(L\right)}m_{\lambda}\left[X\right].
\end{align*}

By \cite[Corollary~3.9]{BAHM26}, for any partition
$\lambda\vdash k$ and any rearrangement $\beta$ of its
parts, we have
\[
    \sum_{L\in\LD(\mm^\bullet;\beta)}q^{\asc(L)}
    =
    \sum_{L\in\LD(\mm^\bullet;\lambda)}q^{\asc(L)}.
\]
Consequently, $\Xi_{\mm^\bullet}[X;q]$ also admits the
quasisymmetric expansion
\[
    \Xi_{\mm^\bullet}[X;q]
    =
    \sum_{\beta\models k}
    \sum_{L\in\LD(\mm^\bullet;\beta)}
    q^{\asc(L)}M_\beta[X].
\]

\begin{rmk}\label{rem:cyclic-chromatic-geometry}
The cyclic chromatic symmetric function has a geometric
interpretation in terms of cyclic Hessenberg varieties.
Let $Z_{\mm^\bullet}$ be the regular semisimple cyclic
Hessenberg variety associated with $\mm^\bullet$ as
constructed in \cite[Section~2]{BAHM26}.
Its cohomology is concentrated in even degrees and
carries an $S_k$-action arising from monodromy.
By \cite[Theorem~3.29]{BAHM26},
\[
    \Xi_{\mm^\bullet}[X;q]
    =
    \sum_{j\geq0}q^j
    \operatorname{Frob}_{S_k}
    \bigl(H^{2j}(Z_{\mm^\bullet})\bigr).
\]
Moreover, let $T\cong(\mathbb C^\times)^k$ be its full
centralizer torus. By \cite[Corollary~3.30]{BAHM26},
\[
    \Xi_{\mm^\bullet}
    \left[\frac{X}{1-q};q\right]
    =
    \sum_{j\geq0}q^j
    \operatorname{Frob}_{S_k}
    \bigl(H_T^{2j}(Z_{\mm^\bullet})\bigr).
\]
Here monodromy permutes the equivariant parameters in
$H_T^*(\mathrm{pt})$, and $q$ records half the
cohomological degree. Thus $\Xi_{\mm^\bullet}[X;q]$
is Schur-positive, and its plethystic transform is
Schur-positive as a formal $q$-series.
\end{rmk}

\begin{lem}\label{lem:omega-rational-chromatic}
For $\mm^\bullet\in\Hessen\left(k,n\right)$, we have
\[
    \omega\left(\Xi_{\mm^\bullet}\left[X;q\right]\right)
    =
    \widetilde{\Xi}_{\mm^\bullet}\left[X;q\right].
\]
\end{lem}

\begin{proof}
For a loop decomposition
$L=\left(L_1,\ldots,L_k\right)\in\LD\left(\mm^\bullet\right)$, define
\[
    G\left(L\right)\coloneqq\left\{i\in\left[k-1\right]: L_i\prec L_{i+1}\right\}.
\]
By the definitions of $\beta$-suitable and $\beta$-colorable, we have
\begin{align*}
     L\in \LD\left(\mm^\bullet;\beta\right)
    &\quad\Longleftrightarrow\quad
    G\left(L\right)\subseteq \Set\left(\beta\right),\\
     L\in \overline{\LD}\left(\mm^\bullet;\beta\right)
    &\quad\Longleftrightarrow\quad
    G\left(L\right)^c\subseteq \Set\left(\beta\right),
\end{align*}
where $G\left(L\right)^c=\left[k-1\right]\setminus G\left(L\right)$.
Recall that, in degree $k$, the fundamental quasisymmetric
function satisfies
\[
    F_S\left[X\right]
    =
    \sum_{\substack{\beta\models k\\ S\subseteq \Set\left(\beta\right)}}
    M_\beta\left[X\right].
\]
We therefore get
\[
\begin{aligned}
    \Xi_{\mm^\bullet}\left[X;q\right]
    &=
    \sum_{\beta\models k}
    \sum_{L\in \LD\left(\mm^\bullet;\beta\right)}
        q^{\asc\left(L\right)}M_\beta\left[X\right]  
    =
    \sum_{L\in\LD\left(\mm^\bullet\right)}
        q^{\asc\left(L\right)}
        \sum_{\substack{\beta\models k\\G\left(L\right)\subseteq\Set\left(\beta\right)}}
        M_\beta\left[X\right] =
    \sum_{L\in\LD\left(\mm^\bullet\right)}
        q^{\asc\left(L\right)}F_{G\left(L\right)}\left[X\right].
\end{aligned}
\]
Applying $\omega$ using \eqref{eq: involution on F} gives
\[
\begin{aligned}
    \omega\left(\Xi_{\mm^\bullet}\left[X;q\right]\right)
    &=
    \sum_{L\in\LD\left(\mm^\bullet\right)}
        q^{\asc\left(L\right)}F_{G\left(L\right)^c}\left[X\right]  =
    \sum_{L\in\LD\left(\mm^\bullet\right)}
        q^{\asc\left(L\right)}
        \sum_{\substack{\beta\models k\\G\left(L\right)^c\subseteq\Set\left(\beta\right)}}
        M_\beta\left[X\right] \\
    &=
    \sum_{\beta\models k}
    \sum_{L\in\overline{\LD}\left(\mm^\bullet;\beta\right)}
        q^{\asc\left(L\right)}M_\beta\left[X\right]=\widetilde{\Xi}_{\mm^\bullet}\left[X;q\right].
\end{aligned}
\]
\end{proof}

\begin{example}
Let $\mmb=\left(\left(1,3,3\right),\left(2,3,3\right)\right)\in \Hessen\left(3,2\right)$. The corresponding directed graph $G\left(\mmb\right)$ is shown in Figure~\ref{fig: Graph mmb}, 
where each vertex $\left(i,3\right)$ is identified with $\left(i,1\right)$. 

\begin{figure}
    \centering
\begin{tikzpicture}[
    x=3.2cm, y=2.6cm,
    >={Latex[length=2.2mm]},
    line width=0.8pt,
    every node/.style={font=\normalsize},
    dot/.style={circle, fill=black, inner sep=3.5pt}
]

\node[dot] (11) at (0,0) {};
\node[dot] (12) at (0,1) {};
\node[dot] (13) at (0,2) {};

\node[dot] (21) at (1,0) {};
\node[dot] (22) at (1,1) {};
\node[dot] (23) at (1,2) {};

\node[dot] (31) at (2,0) {};
\node[dot] (32) at (2,1) {};
\node[dot] (33) at (2,2) {};

\node[below=9pt] at (11) {$\left(1,1\right)$};
\node[below=9pt] at (21) {$\left(2,1\right)$};
\node[below=9pt] at (31) {$\left(3,1\right)$};

\node[below=9pt] at (12) {$\left(1,2\right)$};
\node[below=9pt] at (22) {$\left(2,2\right)$};
\node[below=9pt] at (32) {$\left(3,2\right)$};

\node[above=9pt] at (13) {$\left(1,3\right)$};
\node[above=9pt] at (23) {$\left(2,3\right)$};
\node[above=9pt] at (33) {$\left(3,3\right)$};

\draw[->] (11) -- (12);
\draw[->] (12) -- (13);

\draw[->] (21) -- (22);
\draw[->] (22) -- (23);

\draw[->] (31) -- (32);
\draw[->] (32) -- (33);

\draw[->] (11) -- (22);
\draw[->] (11) -- (32);

\draw[->] (21) -- (32);
\draw[->] (31) -- (22);

\draw[->] (12) -- (23);
\draw[->] (12) -- (33);

\draw[->] (22) -- (13);
\draw[->] (22) -- (33);

\draw[->] (32) -- (23);

\end{tikzpicture}

    \caption{Directed graph $G\left(\mmb\right)$ for $\mmb=\left(\left(1,3,3\right),\left(2,3,3\right)\right)\in \Hessen\left(3,2\right)$.}
    \label{fig: Graph mmb}
\end{figure}
There are two possible (unordered) loop decompositions $\left\{L_1,L_2,L_3\right\}$ and $\left\{L'_1,L'_2,L'_3\right\}$, where
\begin{align*}
    L_1=\left\{\left(1,1\right),\left(1,2\right)\right\}, L_2=\left\{\left(2,1\right),\left(2,2\right)\right\}, L_3=\left\{\left(3,1\right),\left(3,2\right)\right\},\\
    L'_1=\left\{\left(1,1\right),\left(1,2\right)\right\}, L'_2=\left\{\left(2,1\right),\left(3,2\right)\right\}, L'_3=\left\{\left(3,1\right),\left(2,2\right)\right\}.
    \end{align*}
    Note that $L_i\prec L_j$ only when $(i,j)=(3,1)$, whereas $L'_i\prec L'_j$ only when $(i,j)=(2,1)$. Both $\LD\left(\mmb;\left(1,1,1\right)\right)$ and $\overline{\LD}\left(\mmb;\left(1,1,1\right)\right)$ consist of the 12 orderings of $\left(L_1,L_2,L_3\right)$ and $\left(L'_1,L'_2,L'_3\right)$. Consequently we have
        \begin{align*}
        \LD\left(\mmb;\left(2,1\right)\right)&=\LD\left(\mmb;\left(1,1,1\right)\right)\setminus\left\{\left(L_3,L_1,L_2\right),
        \left(L_2',L_1',L_3'\right)\right\},\\
        \LD\left(\mmb;\left(1,2\right)\right)&=\LD\left(\mmb;\left(1,1,1\right)\right)\setminus\left\{\left(L_2,L_3,L_1\right), \left(L'_3,L_2',L_1'\right)\right\}, \text{ and}
        \\\LD\left(\mmb;\left(3\right)\right)&=\LD\left(\mmb;\left(1,1,1\right)\right)\setminus\left\{\left(L_3,L_1,L_2\right),\left(L_2,L_3,L_1\right),\left(L_2',L_1',L_3'\right),\left(L_3',L_2',L_1'\right)\right\}.
    \end{align*}
    We also have 
    \begin{align*}
        \overline{\LD}\left(\mmb;\left(2,1\right)\right)=\left\{\left(L_3,L_1,L_2\right), 
        \left(L_2',L_1',L_3'\right)\right\},\quad 
        \overline{\LD}\left(\mmb;\left(1,2\right)\right)=\left\{\left(L_2,L_3,L_1\right), \left(L'_3,L_2',L_1'\right)\right\},
    \end{align*}
    and $\overline{\LD}\left(\mmb;\left(3\right)\right)$ is empty. Computing the ascent statistic gives
    \begin{align*}
        \Xi_{\mmb}\left[X;q\right]&=\left(q^3+3q^2+3q+1\right)m_{3}+\left(q^3+4q^2+4q+1\right)m_{\left(2,1\right)}+\left(q^3+5q^2+5q+1\right)m_{\left(1,1,1\right)},\\
        \widetilde{\Xi}_{\mmb}\left[X;q\right]&=\left(q^2+q\right)m_{\left(2,1\right)}+\left(q^3+5q^2+5q+1\right)m_{\left(1,1,1\right)}.
    \end{align*}
\end{example}
 
\subsection{Combinatorics of an affine Springer fiber}

In the normalization of Section~\ref{sec:affine-springer},
the formula of Boixeda Alvarez--Hameister--Ma
\cite[Corollary~5.25]{BAHM26} takes the following form.
\begin{equation*}
\mathcal{F}_{km,kn}\left[X,Y;q,t\right]
=
\sum_{\lambda\vdash kn}
m_\lambda\left[X\right]
\sum_{\substack{\pi^\bullet\in\overline{\cPF}_{m,n}^k\left(\lambda\right)\\
    \pi^\bullet\text{ dominant}
}}
t^{\area\left(\pi^\bullet\right)}
q^{\overline{\stat}\left(\pi^\bullet\right)}
\frac{
    \Xi_{\mmb\left(\pi^\bullet\right)}
    \!\left[\dfrac{Y}{1-q};q\right]
}{
    \left[\alpha'\left(\pi^\bullet\right)\right]_q!
}.
\end{equation*}
We now explain the notation appearing in this formula.
\begin{itemize}
    \item For a composition \(\alpha\models kn\), we write \(\overline{\cPF}_{m,n}^k(\alpha)\) for the set of tuples $\pib\in\overline{\cPF}_{m,n}^{k}$
    satisfying
    $X_{\pib}=\prod_{i=1}^{\ell\left(\alpha\right)}x_i^{\alpha_i}$.
    \item   We say that $\pi^{\bullet}\in \overline{\cPF}^{k}_{m,n}\left(\alpha\right)$ is \emph{dominant} if
    \begin{equation*}
        \pi^{\left(1\right)}_{\ell}\leq\pi^{\left(2\right)}_{\ell}\leq \cdots \leq \pi^{\left(k\right)}_{\ell}
    \end{equation*}
    for each $1\leq \ell\leq n$. Here, the order is the total order in \eqref{eq:labeled-step-order}.
\item For a dominant $\pi^\bullet\in \overline{\cPF}^{k}_{m,n}\left(\alpha\right)$,
we associate $\mmb\left(\pi^\bullet\right)\in\Hessen\left(k,n\right)$ as follows. For each $1\leq \ell\leq n$ and $1 \leq j\leq k$, $\mm^{\left(\ell\right)}_{j}$ is defined to be the largest integer  $i\in[k]$ such that $\pi^{\left(k-i+1\right)}_{\ell}$ can be weakly extended with $\pi^{\left(k-j+1\right)}_{\ell+1}$. Dominance implies that each $\mm^{(\ell)}$ is weakly increasing.
Moreover, $i=j$ is admissible by the weak cyclic condition, so
$j\leq\mm_j^{(\ell)}\leq k$.
Thus $\mm^{(\ell)}$ is a Hessenberg function of size $k$.

\item We define a partition $\alpha'\left(\pi^\bullet\right)\vdash kn$ by recording the
multiplicities of the labeled north-step types in $\pi^\bullet$
(with the sheet index forgotten) and then rearranging
these multiplicities in weakly decreasing order.
\end{itemize}

\begin{rmk}
To match the dominant convention of~\cite{BAHM26},
we reverse the order of the $k$ components before assigning
contents.
The $kn$-generators of the corresponding invariant set
are $\cont(N)+1$, with labels inherited from the north steps. This correspondence preserves area and identifies the
parabolic automorphism factor with
$[\alpha'(\pi^\bullet)]_q!$.

For corresponding dominant representatives, the statistics satisfy
\[
    \overline{\stat}(\pi^\bullet)
    =
    \delta-\operatorname{codinv}(\Delta,f)-r,
    \qquad
    r=
    D_{\mmb(\pi^\bullet)}
    -\sum_j\binom{\alpha'_j(\pi^\bullet)}{2}.
\]
Here $r$ is the dimension of the associated parabolic cyclic
Hessenberg variety.
Together with Poincar\'e duality, this translates
\cite[Corollary~5.25]{BAHM26} into
\eqref{eq: affine Springer fiber formula}
with our grading conventions.
\end{rmk}

Let
\[
    \skel\left(m,n,k;\alpha\right)
    \coloneqq
    \left\{z\left(\pi^\bullet\right):
      \pi^\bullet\in\overline{\cPF}_{m,n}^{k}\left(\alpha\right)\right\}
\]
be the set of skeletons of elements of
\(\overline{\cPF}_{m,n}^{k}(\alpha)\).
Each skeleton $z\in\skel(m,n,k;\alpha)$ has a unique
dominant representative
$\pi^\bullet\in\overline{\cPF}_{m,n}^k(\alpha)$.
We therefore identify skeletons with their dominant representatives.
\begin{example}
Let $\pib=\left(\pi^{\left(1\right)},\pi^{\left(2\right)},\pi^{\left(3\right)}\right)\in\overline{\cPF}^{3}_{3,2}$ be the tuple shown in Figure~\ref{fig: dominant}. Note that $\pi^{\left(i\right)}_{1}$ can be weakly extended with  $\pi^{\left(j\right)}_{2}$ for all \(i\) and \(j\). Meanwhile, $\pi^{\left(i\right)}_{2}$ can be weakly extended with  $\pi^{\left(j\right)}_{3}$ except when $i=1$ and $j=3$. We conclude
\begin{equation*}
    \mmb\left(\pib\right)=\left(\left(3,3,3\right),\left(2,3,3\right)\right).
\end{equation*}
We have $\alpha'\left(\pib\right)=\left(2,1,1,1,1\right)$ since $\pi^{\left(2\right)}_2=\pi^{\left(3\right)}_2$ and there are no other repeated labeled north-step types.
    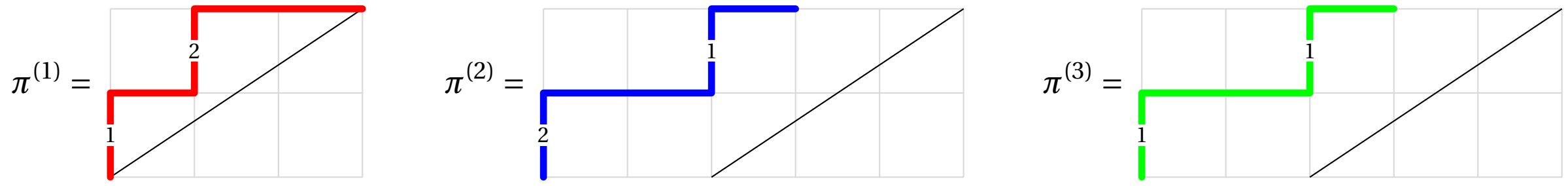
\begin{figure}[ht]
\centering

$\displaystyle \pi^{\left(1\right)}=$
\begin{tikzpicture}[scale=0.85, baseline=(current bounding box.center)]
    
    \draw[help lines, gray!30, thin] (0,0) grid (3,2);
    \draw[black, thin] (0,0) -- (3,2);

    \draw[red, line width=2pt, line cap=round, line join=round]
        (0,0) -- (0,1) -- (1,1)
        -- (1,2) -- (2,2) -- (3,2);

    \node[fill=white, inner sep=1pt, text=black, font=\scriptsize]
        at (0,0.5) {$1$};
    \node[fill=white, inner sep=1pt, text=black, font=\scriptsize]
        at (1,1.5) {$2$};
\end{tikzpicture}
\qquad
$\displaystyle \pi^{\left(2\right)}=$
\begin{tikzpicture}[scale=0.85, baseline=(current bounding box.center)]

    \draw[help lines, gray!30, thin] (-2,0) grid (3,2);
    \draw[black, thin] (0,0) -- (3,2);

    \draw[blue, line width=2pt, line cap=round, line join=round]
        (-2,0) -- (-2,1) -- (0,1)
        -- (0,2) -- (1,2);

    \node[fill=white, inner sep=1pt, text=black, font=\scriptsize]
        at (-2,0.5) {$2$};
    \node[fill=white, inner sep=1pt, text=black, font=\scriptsize]
        at (0,1.5) {$1$};
\end{tikzpicture}
\qquad
$\displaystyle \pi^{\left(3\right)}=$
\begin{tikzpicture}[scale=0.85, baseline=(current bounding box.center)]

    \draw[help lines, gray!30, thin] (-2,0) grid (3,2);
    \draw[black, thin] (0,0) -- (3,2);

    \draw[green, line width=2pt, line cap=round, line join=round]
        (-2,0) -- (-2,1) -- (0,1)
        -- (0,2) -- (1,2);

    \node[fill=white, inner sep=1pt, text=black, font=\scriptsize]
        at (-2,0.5) {$1$};
    \node[fill=white, inner sep=1pt, text=black, font=\scriptsize]
        at (0,1.5) {$1$};
\end{tikzpicture}

\caption{A dominant tuple $\pib=\left(\pi^{\left(1\right)},\pi^{\left(2\right)},\pi^{\left(3\right)}\right)\in\overline{\cPF}^{3}_{3,2}$.}
\label{fig: dominant}
\end{figure}
\end{example}

Fix a dominant $\pi^\bullet\in\overline{\cPF}_{m,n}^{k}(\alpha)$
and a composition $\beta\models k$.
For a skeleton $\z$, write
\[
    \overline{\CE}_\beta\left(\z\right)
    \coloneqq
    \left\{
        \tau^\bullet\in\overline{\CE}_\beta:
        z\left(\tau^\bullet\right)=\z
    \right\}.
\]
Set $\mmb=\mmb\left(\pi^\bullet\right)$, and let
\[
    V=V\left(G\left(\mmb\right)\right)
    =
    \left\{(i,\ell):1\leq i\leq k,\ 1\leq\ell\leq n\right\}.
\]
For simplicity, write
\[
    \gamma\left(i,\ell\right)\coloneqq\pi^{\left(k-i+1\right)}_\ell
\]
for the labeled north step associated with the vertex $\left(i,\ell\right)$.

Partition $V$ into equivalence classes $I_1,\ldots,I_c$,
where $\left(i,\ell\right)$ and $\left(i',\ell'\right)$ belong to the same class
if and only if $\gamma\left(i,\ell\right)=\gamma\left(i',\ell'\right)$.
Fix an ordering of these classes such that
\[
    \alpha'\left(\pi^\bullet\right)=\left(\left|I_1\right|,\ldots,\left|I_c\right|\right).
\]
We write
\[
    \mathfrak S_{\alpha'\left(\pi^\bullet\right)}
    =
    \mathfrak S_{\left|I_1\right|}\times\cdots\times\mathfrak S_{\left|I_c\right|}
\]
for the corresponding Young subgroup.

Given $L=\left(L_1,\ldots,L_k\right)\in\LD\left(\mmb;\beta\right)$, define
$\tau^\bullet=\left(\tau^{\left(1\right)},\ldots,\tau^{\left(k\right)}\right)$ by
\[
    \tau^{\left(i\right)}_\ell
    =
    \gamma\left(L_i\left(\ell\right),\ell\right)
    =
    \pi^{\left(k-L_i\left(\ell\right)+1\right)}_\ell
    \qquad \left(1\leq i\leq k,\ 1\leq\ell\leq n\right).
\]
Since each $L_i$ is a loop, we have
$\tau^{\left(i\right)}\in\overline{\cPF}_{m,n}$.
Moreover, $\tau^\bullet$ has the same skeleton as $\pi^\bullet$,
and the $\beta$-suitability of $L$ implies that
\[
    \tau^\bullet\in
    \overline{\CE}_\beta\left(z\left(\pi^\bullet\right)\right).
\]

By dominance, each class $I_j$ is of the form
\[
    I_j=\left\{\left(i,\ell\right),\left(i+1,\ell\right),\ldots,\left(i+p,\ell\right)\right\}.
\]
For each $0\leq s\leq p$, let $d_s$ be the unique index such that
\[
    \left(i+s,\ell\right)\in L_{d_s}.
\]
Define $w^{\left(j\right)}\in\mathfrak S_{\left|I_j\right|}$ by standardizing
$\left(d_p,d_{p-1},\ldots,d_0\right)$, and set
\[
    F\left(L\right)\coloneqq\left(\tau^\bullet,w^\bullet\right),
    \qquad
    w^\bullet=\left(w^{\left(1\right)},\ldots,w^{\left(c\right)}\right).
\]

\begin{proposition}\label{prop: L to pib}
The map
\[
    F:
    \LD\left(\mmb\left(\pi^\bullet\right);\beta\right)
    \longrightarrow
    \overline{\CE}_\beta\left(z\left(\pi^\bullet\right)\right)
    \times\mathfrak S_{\alpha'\left(\pi^\bullet\right)}
\]
defined above is a bijection such that
\[
    \asc\left(L\right)
    =
    \overline{\pos}\left(\tau^\bullet\right)
    +\sum_{j=1}^{c}\inv\left(w^{\left(j\right)}\right),
\]
whenever $F\left(L\right)=\left(\tau^\bullet,w^\bullet\right)$.
\end{proposition}

\begin{proof}
We first verify the statistic identity.
Suppose $\left(\ell,i,j\right)$ is an ascent of $L$, and put
$a=L_i\left(\ell\right)$ and $b=L_j\left(\ell\right)$.
Then $i<j$ and
\[
    a<b\leq\mm^{\left(\ell\right)}_{L_i\left(\ell+1\right)}.
\]
\begin{enumerate}
    \item If $\gamma\left(a,\ell\right)\ne\gamma\left(b,\ell\right)$, then
    $\left(\ell,i,j\right)$ is a weak positive crossing of $\tau^\bullet$.

    \item If $\gamma\left(a,\ell\right)=\gamma\left(b,\ell\right)$, let $I_g$
    be the class containing $\left(a,\ell\right)$ and $\left(b,\ell\right)$.
    The occurrences assigned to loops $i$ and $j$
    contribute one inversion to $w^{\left(g\right)}$.
\end{enumerate}
Conversely, every weak positive crossing of $\tau^\bullet$ and
every inversion of some $w^{\left(g\right)}$ arises from a unique
ascent in the corresponding case.
Thus
\[
    \asc\left(L\right)
    =
    \overline{\pos}\left(\tau^\bullet\right)+\sum_{j=1}^{c}\inv\left(w^{\left(j\right)}\right).
\]

To prove bijectivity, let
\[
    \left(\tau^\bullet,w^\bullet\right)
    \in
    \overline{\CE}_\beta\left(z\left(\pi^\bullet\right)\right)
    \times\mathfrak S_{\alpha'\left(\pi^\bullet\right)}.
\]
Within each equivalence class $I_j$, the tuple \(\tau^\bullet\) determines the set of loop indices to which the corresponding labeled north-step type is assigned, and $w^{\left(j\right)}$ uniquely determines their
assignment to the vertices of $I_j$ in the reversed order
used above.
The weak cyclic condition on $\tau^\bullet$ ensures that
the resulting $L_i$ are loops of $G\left(\mmb\right)$, and the condition
$\tau^\bullet\in\overline{\CE}_\beta$ is equivalent to the
$\beta$-suitability of $L$.
This reconstruction is inverse to $F$, so $F$ is a bijection.
\end{proof}

\begin{proof}[Proof of Theorem~\ref{thm: main}]
By Proposition~\ref{prop: omega Omega}, it suffices to prove 
\begin{equation}
    \sum_{\mu\vdash k}
    \left(\mathfrak e^{\mathrm{w}}_\mu\cdot 1\right)
    m_\mu\!\left[\frac{Y}{1-q}\right]=\sum_{\lambda \vdash kn}
m_{\lambda}\left[X\right]\sum_{\substack{\pi^{\bullet}\in \overline{\cPF}_{m,n}^{k}\left(\lambda\right)\\\pi^{\bullet}\text{ dominant}}}
t^{\area\left(\pi^{\bullet}\right)}
q^{\overline{\stat}\left(\pi^{\bullet}\right)}
\frac{
\Xi_{\mmb\left(\pi^{\bullet}\right)}\!\left[\dfrac{Y}{1-q};q\right]
}{
\left[\alpha'\left(\pi^{\bullet}\right)\right]_{q}!}.
\end{equation}
Fix $\lambda\vdash kn$. Taking the coefficient of
$m_\lambda\left[X\right]$ reduces the desired identity to
\begin{equation*}
    \sum_{\mu\vdash k}
    \left(
        \sum_{\z\in\skel\left(m,n,k;\lambda\right)}
        \sum_{\substack{
            \tau^\bullet\in\overline{\CE}_\mu\\
            z\left(\tau^\bullet\right)=\z
        }}
        t^{\area\left(\tau^\bullet\right)}
        q^{\overline{\stat}\left(\tau^\bullet\right)}
    \right)
    m_\mu\!\left[\frac{Y}{1-q}\right]
    =
    \sum_{\substack{
        \pi^\bullet\in\overline{\cPF}_{m,n}^k\left(\lambda\right)\\
        \pi^\bullet\text{ dominant}
    }}
    t^{\area\left(\pi^\bullet\right)}
    q^{\overline{\stat}\left(\pi^\bullet\right)}
    \frac{
        \Xi_{\mmb\left(\pi^\bullet\right)}
        \!\left[\dfrac{Y}{1-q};q\right]
    }{
        \left[\alpha'\left(\pi^\bullet\right)\right]_q!
    }.
\end{equation*}
Comparing coefficients in the basis
$\left\{m_\mu\left[Y/(1-q)\right]:\mu\vdash k\right\}$,
it remains to show that, for every $\mu\vdash k$,
\begin{equation*}
    \sum_{\z\in\skel\left(m,n,k;\lambda\right)}
    \sum_{\substack{
        \tau^\bullet\in\overline{\CE}_\mu\\
        z\left(\tau^\bullet\right)=\z
    }}
    t^{\area\left(\tau^\bullet\right)}
    q^{\overline{\stat}\left(\tau^\bullet\right)}
    =
    \sum_{\substack{
        \pi^\bullet\in\overline{\cPF}_{m,n}^k\left(\lambda\right)\\
        \pi^\bullet\text{ dominant}
    }}
    t^{\area\left(\pi^\bullet\right)}
    q^{\overline{\stat}\left(\pi^\bullet\right)}
    \frac{
        \displaystyle
        \sum_{L\in\LD\left(\mmb\left(\pi^\bullet\right);\mu\right)}
        q^{\asc\left(L\right)}
    }{
        \left[\alpha'\left(\pi^\bullet\right)\right]_q!
    }.
\end{equation*}
Fix a dominant \(\pi^\bullet\) and set \(\z=z\left(\pi^\bullet\right)\). By
Proposition~\ref{prop: L to pib},
\[
\sum_{L\in\LD\left(\mmb\left(\pi^\bullet\right);\mu\right)}q^{\asc\left(L\right)}
=
\left[\alpha'\left(\pi^\bullet\right)\right]_q!
\sum_{\substack{
    \tau^\bullet\in\overline{\CE}_\mu\\
    z\left(\tau^\bullet\right)=\z
}}
q^{\overline{\pos}\left(\tau^\bullet\right)}.
\]
Since \(\pi^\bullet\) is dominant, \(\overline{\pos}\left(\pi^\bullet\right)=0\). Moreover,
Lemma~\ref{lem: pos crossing} implies
\[
    \overline{\stat}\left(\tau^\bullet\right)
    =
    \overline{\stat}\left(\pi^\bullet\right)
    +\overline{\pos}\left(\tau^\bullet\right)
\]
for all \(\tau^\bullet\) with \(z\left(\tau^\bullet\right)=\z\). Since such
\(\tau^\bullet\) have the same area as \(\pi^\bullet\), the desired identity
follows.
\end{proof}

\section{Cyclic LLT polynomials and Theorem~\ref{thm: main B}}\label{Sec: rational uLLT}

\subsection{Cyclic LLT polynomials}

The plethystic transforms appearing in Theorem~\ref{thm: main} lead to
the following cyclic analogues of unicellular LLT polynomials.

\begin{definition}
For $\mmb\in\Hessen(k,n)$, define the
\emph{cyclic LLT polynomial} by
\[
    \cLLT_{\mmb}[X;q]
    \coloneqq
    (1-q)^k\,
    \Xi_{\mmb}\!\left[\frac{X}{1-q};q\right].
\]
\end{definition}

When $n=1$, the tuple $\mmb\in\Hessen(k,1)$ consists of
a single Hessenberg function $\mm^{(1)}$, and
$\cLLT_{\mmb}[X;q]$ recovers the usual unicellular LLT
polynomial associated with $\mm^{(1)}$.

By Remark~\ref{rem:cyclic-chromatic-geometry}, the quotient
$\cLLT_{\mmb}[X;q]/(1-q)^k$ is Schur-positive as a formal
$q$-series. Computations suggest the following stronger
statement.

\begin{conj}
For every $\mmb\in\Hessen(k,n)$, the cyclic LLT polynomial
$\cLLT_{\mmb}[X;q]$ is Schur-positive, i.e. its Schur coefficients lie in $\mathbb{Z}_{\ge0}[q]$.
\end{conj}

Finding a positive combinatorial formula for the Schur coefficients of
\(\cLLT_{\mmb}\left[X;q\right]\) appears to be difficult: even in the unicellular LLT case, a general positive combinatorial formula for the Schur coefficients remains open. In Theorem~\ref{thm: rLLt monomial} below, we prove a monomial expansion of
\(\cLLT_{\mmb}\left[X;q\right]\), generalizing the usual monomial formula for unicellular LLT polynomials. Before stating the formula, we introduce the necessary terminology.

For a composition $\beta\models k$, set
$r_i\coloneqq1+\sum_{j=1}^{i-1}\beta_j$ for each
$1\le i\le\ell\left(\beta\right)$.
Consider a loop decomposition
$L=\left(L_1,L_2,\dots,L_k\right)\in\LD\left(\mm^\bullet\right)$.
\begin{itemize}
    \item We say that $L$ is \emph{$\beta$-sorted} if
    \begin{equation*}
        L_{r_i}\left(1\right)>L_{r_i+1}\left(1\right)>\cdots>L_{r_i+\beta_i-1}\left(1\right)
    \end{equation*}
    for each $1\le i\le\ell\left(\beta\right)$.

    \item We say that $L$ is \emph{reverse $\beta$-sorted} if
    \begin{equation*}
        L_{r_i}\left(1\right)<L_{r_i+1}\left(1\right)<\cdots<L_{r_i+\beta_i-1}\left(1\right)
    \end{equation*}
    for each $1\le i\le\ell\left(\beta\right)$.
\end{itemize}

 We denote the sets of $\beta$-sorted and reverse
$\beta$-sorted loop decompositions by
$\sLD(\mm^\bullet;\beta)$ and
$\overline{\sLD}(\mm^\bullet;\beta)$, respectively.

\begin{thm}\label{thm: rLLt monomial} Let $k,n$ be positive integers and let
$\mmb\in\Hessen(k,n)$. Then we have
    \begin{align}
        \cLLT_{\mm^{\bullet}}\left[X;q\right]=\left(1-q\right)^k \Xi_{\mmb}\left[\frac{X}{1-q};q\right]&=\sum_{\lambda \vdash k}\sum_{L\in \sLD\left(\mm^{\bullet};\lambda\right)}q^{\asc\left(L\right)}m_{\lambda}\label{eq: llt},\\
        \left(1-q\right)^k \widetilde{\Xi}_{\mmb}\left[\frac{X}{1-q};q\right]&=\sum_{\lambda \vdash k}\sum_{L\in \overline{\sLD}\left(\mm^{\bullet};\lambda\right)}q^{\asc\left(L\right)}m_{\lambda}.\label{eq: dual llt}
    \end{align}
\end{thm}

Note that when $n=1$, Theorem~\ref{thm: rLLt monomial} reduces to the usual monomial expansion for unicellular LLT polynomials. 

To prove Theorem~\ref{thm: rLLt monomial} from the monomial
expansions of $\Xi_{\mmb}[X;q]$ and $\widetilde{\Xi}_{\mmb}[X;q]$,
we first compute
$[m_\mu]m_\lambda[X/(1-q)]$.
By the coefficient identity \eqref{eq: e m}, this equals
$[e_\lambda]e_\mu[X/(1-q)]$.

For a composition \(\alpha=\left(\alpha_1,\dots,\alpha_\ell\right)\models n\), we define
\begin{equation}\label{eq: S alpha define}
    S\left(\alpha\right)
    =
    q^{\sum_{i\in\Set\left(\alpha\right)}i}
    \prod_{i\in \Set\left(\alpha\right)^c}\left(1-q^i\right).
\end{equation}

\begin{proposition}\label{prop: en to e}
For \(n\ge 1\), we have
\[
    e_n\!\left[\frac{X}{1-q}\right]
    =
    \sum_{\alpha\models n}
\frac{S\left(\alpha\right)}{\left(1-q\right)^n\left[n\right]_q!}e_\alpha.
\]
\end{proposition}

\begin{proof} 
By considering \(X/\left(1-q\right)=X+qX+q^2X+\cdots\), we obtain 
\[
    e_n\!\left[\frac{X}{1-q}\right]
    =
    \sum_{\alpha\models n} \left(\sum_{r_1>\cdots>r_\ell\ge 0}
    q^{r_1\alpha_1+\cdots+r_\ell\alpha_\ell}\right)
    e_\alpha.
\]
Set
\(
    W_\alpha = \sum_{r_1>\cdots>r_\ell\ge 0}
    q^{r_1\alpha_1+\cdots+r_\ell\alpha_\ell}
\)
and \(a_j=\alpha_1+\cdots+\alpha_j\).  If
\(s_j=r_j-r_{j+1}\) for \(1\le j<\ell\) and \(s_\ell=r_\ell\), then
\(s_1,\dots,s_{\ell-1}\ge 1\), \(s_\ell\ge 0\), and
\[
r_1\alpha_1+\cdots+r_\ell\alpha_\ell
=
s_1a_1+\cdots+s_\ell a_\ell .
\]
Since $a_\ell=n$ and
$\Set\left(\alpha\right)=\left\{a_1,\ldots,a_{\ell-1}\right\}$, we obtain
\[
\begin{aligned}
    W_\alpha
    &=
    \left(\sum_{s_\ell\geq0}q^{s_\ell n}\right)
    \prod_{j=1}^{\ell-1}
    \left(\sum_{s_j\geq1}q^{s_j a_j}\right)=
    \frac{1}{1-q^n}
    \prod_{j=1}^{\ell-1}\frac{q^{a_j}}{1-q^{a_j}}=
    \frac{S\left(\alpha\right)}{\prod_{i=1}^n\left(1-q^i\right)}.
\end{aligned}
\]
Therefore,
\[
e_n\!\left[\frac{X}{1-q}\right]
=
\sum_{\alpha\models n}
\frac{S\left(\alpha\right)}{\left(1-q\right)^n\left[n\right]_q!}e_\alpha .
\]
\end{proof}

Combining Proposition~\ref{prop: en to e} with the coefficient
identity \eqref{eq: e m} gives the following coefficient formula.

\begin{corollary}\label{cor: m plethysm}
Let $\lambda,\mu\vdash n$, and let
$\beta=(\beta_1,\ldots,\beta_\ell)$ be a rearrangement of $\mu$.
Then
\[
    \left[m_\mu\left[X\right]\right]\,m_\lambda\!\left[\frac{X}{1-q}\right]
    =
    \frac{1}{\left(1-q\right)^{n}\left[\beta\right]_q!}
    \sum_{\substack{
        \alpha^{\left(i\right)}\models \beta_i\\
        \left(\alpha^{\left(1\right)},\dots,\alpha^{\left(\ell\right)}\right)\sim\lambda
    }}
    \prod_{i=1}^{\ell} S\left(\alpha^{\left(i\right)}\right),
\]
where \(\left(\alpha^{\left(1\right)},\dots,\alpha^{\left(\ell\right)}\right)\) denotes concatenation.
\end{corollary}

\subsection{Proof of Theorem~\ref{thm: rLLt monomial}}

Throughout this subsection, we fix $\mmb\in\Hessen\left(k,n\right)$.
We first prove the dual identity \eqref{eq: dual llt},
and then deduce \eqref{eq: llt} by applying $\omega$.

For $0\le r\le k$, let $\pLD_r\left(\mm^\bullet\right)$
be the set of ordered $r$-tuples $L=(L_1,\ldots,L_r)$
of pairwise disjoint loops of $G\left(\mm^\bullet\right)$
such that the entire tuple can be extended to a loop decomposition
of $\mm^\bullet$. For \(r=k\), this gives \(\pLD_k\left(\mmb\right)=\LD\left(\mmb\right)\). The statistic
\(\asc\left(L\right)\) is defined by the same ascent condition as for full loop
decompositions, but only among the loops \(L_1,\ldots,L_r\).

Let \(A\subseteq \left[k\right]\) with \(\left|A\right|=r\).
Throughout this subsection, a word on \(A\) means
a permutation of the elements of \(A\).
A \emph{bracketed word} on \(A\) is such a word,
together with a decomposition into consecutive bracketed blocks.
If
\[
    W=\left(w_1,\ldots,w_r\right)
\]
is such a bracketed word, define
\[
    \left\langle W\right\rangle
    \coloneqq
    \sum_{L} q^{\asc\left(L\right)},
\]
where the sum is over all
\(L=\left(L_1,\ldots,L_r\right)\in\pLD_r\left(\mm^\bullet\right)\) satisfying
\(L_i\left(1\right)=w_i\) for every \(i\), and satisfying the following additional
condition: whenever \(\left[w_a,\ldots,w_b\right]\) is a bracketed block of \(W\), we impose
\[
    L_a\prec L_{a+1}\prec\cdots\prec L_b .
\]
A block of length one may be written either with or without brackets. Thus, for
example,
\[
    \left\langle a_1,\left[a_2,a_3\right],a_4\right\rangle
\]
denotes the generating function over all
\(\left(L_1,L_2,L_3,L_4\right)\in\pLD_4\left(\mm^\bullet\right)\) such that
\(L_i\left(1\right)=a_i\) and \(L_2\prec L_3\).

For a composition \(\alpha=\left(\alpha_1,\ldots,\alpha_\ell\right)\models r\), let
\(\mathcal B\left(A;\alpha\right)\) be the set of bracketed words on \(A\) whose bracketed
blocks have lengths \(\alpha_1,\ldots,\alpha_\ell\), from left to right. When
\(A=\left[r\right]\), we write simply \(\mathcal B\left(\alpha\right)\).

We use the following $q$-swapping relation, in which
the bracketed words $U$ and $V$ are fixed.

\begin{lem}\label{lem:q-swapping}
Let $r\geq1$, let $a_1<\cdots<a_{r+1}$ be elements of $\left[k\right]$, and let
$U,V$ be possibly empty bracketed words.
Assume that all letters occurring in $U$, $V$, and
$a_1,\ldots,a_{r+1}$ are pairwise distinct and lie in $\left[k\right]$.
Then
\[
\begin{aligned}
    &\left\langle U,\left[a_r,\ldots,a_1\right],a_{r+1},V\right\rangle
    -q\left\langle U,a_{r+1},\left[a_r,\ldots,a_1\right],V\right\rangle\\
    &\qquad =
    \left(1-q\right)\left\langle U,\left[a_{r+1},a_r,\ldots,a_1\right],V\right\rangle.
\end{aligned}
\]
\end{lem}

This identity is a local form of the straightening relation
from \cite{KO26}.

\begin{example}
    Let $\mm^{\bullet}=\left(\left(1,3,3\right),\left(1,3,3\right)\right)$. Then there are two (unordered) loop decompositions $L=\left\{L_1,L_2,L_3\right\}$ and $L'=\left\{L_1,L_2',L_3'\right\}$, where
    \begin{align*}
        L_1=\{\left(1,1\right),&\left(1,2\right)\}, L_2=\left\{\left(2,1\right),\left(2,2\right)\right\}, L_3=\left\{\left(3,1\right),\left(3,2\right)\right\},\\
        &L_2'=\left\{\left(2,1\right),\left(3,2\right)\right\}, L_3'=\left\{\left(3,1\right),\left(2,2\right)\right\}.
    \end{align*}
    Then we have
    \begin{align*}
      \left\langle\left[2,1\right],3\right\rangle=q^{\asc\left(L_2,L_1,L_3\right)}+q^{\asc\left(L'_2,L_1,L'_3\right)}=q^2+q,\\
      \left\langle 3,\left[2,1\right]\right\rangle=q^{\asc\left(L_3,L_2,L_1\right)}+q^{\asc\left(L_3',L'_2,L_1\right)}=1+q,
    \end{align*}
    and $\left\langle \left[3,2,1\right]\right\rangle=0$. We verify that 
    \begin{equation*}
        \left\langle\left[2,1\right],3\right\rangle-q\left\langle3,\left[2,1\right]\right\rangle=0=\left(1-q\right)\left\langle \left[3,2,1\right]\right\rangle,
    \end{equation*}
    in agreement with Lemma~\ref{lem:q-swapping}.
\end{example}

We next remove brackets by iterating this relation (Lemma~\ref{lem:bracket-deform}).

\begin{definition}
Let \(A=\{a_1<\cdots<a_r\}\subseteq[k]\), where \(r\ge1\). A word
\(v=\left(v_1,\ldots,v_r\right)\) on \(A\) is called \emph{\(V\)-shaped} if there exists
\(d\in\left[r\right]\) such that
\[
    v_1>\cdots>v_d
    \qquad\text{and}\qquad
    v_d<v_{d+1}<\cdots<v_r .
\]
We allow \(d=1\), in which case the word is increasing, and \(d=r\), in which
case the word is decreasing. This integer \(d\) is uniquely determined; we
denote it by \(d\left(v\right)\). Let \(\mathcal V\left(A\right)\) be the set of \(V\)-shaped words on
\(A\).
\end{definition}

The identities in Lemma~\ref{lem:bracket-deform} and
Proposition~\ref{prop:technical-bracket} remain valid
upon replacing every $\langle W\rangle$ by
$\langle U,W,V\rangle$, where $U,V$ are fixed bracketed
words and all letters in $U$, $V$, and $A$ are pairwise
distinct.
The same proofs apply, using Lemma~\ref{lem:q-swapping}
with the fixed words $U,V$ retained throughout.

\begin{lem}\label{lem:bracket-deform}
Let \(A=\{a_1<\cdots<a_r\}\subseteq[k]\), where \(r\ge1\). Then
\[
    \left(1-q\right)^{r-1}
    \left\langle \left[a_r,a_{r-1},\ldots,a_1\right]\right\rangle
    =
    \sum_{v\in\mathcal V\left(A\right)}
        \left(-q\right)^{d\left(v\right)-1}\left\langle v\right\rangle .
\]
\end{lem}

\begin{proof}
We proceed by induction on \(r\). The case \(r=1\) is immediate. For \(r\ge2\),
Lemma~\ref{lem:q-swapping} gives
\[
\begin{aligned}
    \left(1-q\right)
    \left\langle \left[a_r,a_{r-1},\ldots,a_1\right]\right\rangle
    &=
    \left\langle \left[a_{r-1},\ldots,a_1\right],a_r\right\rangle  \\
    &\qquad
    -q\left\langle a_r,\left[a_{r-1},\ldots,a_1\right]\right\rangle .
\end{aligned}
\]
Multiplying by \(\left(1-q\right)^{r-2}\) and applying the induction hypothesis in the two
fixed contexts gives
\[
\begin{aligned}
    \left(1-q\right)^{r-1}
    \left\langle \left[a_r,\ldots,a_1\right]\right\rangle
    &=
    \sum_{u\in\mathcal V\left(A\setminus\left\{a_r\right\}\right)}
        \left(-q\right)^{d\left(u\right)-1}\left\langle u,a_r\right\rangle        \\
    &\qquad
    +
    \sum_{u\in\mathcal V\left(A\setminus\left\{a_r\right\}\right)}
        \left(-q\right)^{d\left(u\right)}\left\langle a_r,u\right\rangle .
\end{aligned}
\]
Every $V$-shaped word on $A$ either ends in $a_r$ or begins
with $a_r$, and these two cases are disjoint.
Moreover, if $v=\left(u,a_r\right)$, then $d\left(v\right)=d\left(u\right)$, while if
$v=\left(a_r,u\right)$, then $d\left(v\right)=d\left(u\right)+1$.
Thus the two sums combine to give
\[
    \sum_{v\in\mathcal V\left(A\right)}
    \left(-q\right)^{d\left(v\right)-1}\left\langle v\right\rangle.
\]
\end{proof}

We now prove the main straightening identity needed for the plethystic
substitution. Recall the value $S\left(\alpha\right)$ defined for a composition $\alpha$ in \eqref{eq: S alpha define}.
\begin{proposition}\label{prop:technical-bracket}
Let \(A=\{a_1<\cdots<a_r\}\subseteq[k]\), where \(r\ge1\). Then
\[
    \sum_{\alpha\models r}
        S\left(\alpha\right)
        \sum_{W\in\mathcal B\left(A;\alpha\right)}
            \left\langle W\right\rangle
    =
    \left[r\right]_q!\,\left\langle a_1,a_2,\ldots,a_r\right\rangle .
\]
\end{proposition}

\begin{proof}
Fix \(\alpha=\left(\alpha_1,\ldots,\alpha_\ell\right)\models r\). A word
\(w=\left(w_1,\ldots,w_r\right)\) on \(A\) is called \(\alpha\)-\(V\)-shaped if, after
cutting \(w\) into consecutive blocks of lengths
\(\alpha_1,\ldots,\alpha_\ell\), each block is \(V\)-shaped. Let
\(\mathcal V\left(A;\alpha\right)\) be the set of such words. If
\(w=w^{\left(1\right)}\cdots w^{\left(\ell\right)}\) is the corresponding block decomposition, set
\[
    d_\alpha\left(w\right)
    \coloneqq
    \sum_{j=1}^{\ell} d\left(w^{\left(j\right)}\right).
\]
Applying Lemma~\ref{lem:bracket-deform} to each bracketed block yields
\[
\begin{aligned}
    S\left(\alpha\right)
    \sum_{W\in\mathcal B\left(A;\alpha\right)}
        \left\langle W\right\rangle
    =
    q^{\binom r2}
    \eta\left(\alpha\right)
    \sum_{w\in\mathcal V\left(A;\alpha\right)}
        \left(-q\right)^{d_\alpha\left(w\right)-\ell\left(\alpha\right)}
        \left\langle w\right\rangle ,
\end{aligned}
\]
where
\[
    \eta\left(\alpha\right)
    \coloneqq
    \prod_{i\in \Set\left(\alpha\right)^{c}}
        q^{-i}\left[i\right]_q .
\]
Indeed,
\[
    \frac{S\left(\alpha\right)}{\left(1-q\right)^{r-\ell\left(\alpha\right)}}
    =
    q^{\binom r2}\eta\left(\alpha\right).
\]

We now sum over all $\alpha\models r$ and collect the
contributions indexed by a fixed word $w=(b_1,\ldots,b_r)$ on $A$. Suppose first that \(w\) has a descent,
and let \(i\) be the smallest index such that \(b_i>b_{i+1}\).

If \(i=1\), the compositions contributing to the coefficient of \(\left\langle w\right\rangle\)
are paired by toggling the cut at position \(1\). In each pair, the two local
contributions are
\[
    1
    \qquad\text{and}\qquad
    \left(-q\right)q^{-1}\left[1\right]_q,
\]
whose sum is zero.

Now assume \(i>1\). Since \(i\) is the first descent, we have
\(b_{i-1}<b_i>b_{i+1}\). Hence an \(\alpha\)-\(V\)-shaped decomposition must
have a cut at \(i-1\) or at \(i\). We group the contributing compositions into
triples according to whether the set of cuts among positions \(i-1\) and \(i\) is \(\{i-1,i\}\), \(\{i-1\}\), or \(\{i\}\), with all other cuts fixed. Let these three compositions be
\(\alpha^{\left(1\right)},\alpha^{\left(2\right)},\alpha^{\left(3\right)}\). Then
\[
    \eta\left(\alpha^{\left(2\right)}\right)
    =
    q^{-i}\left[i\right]_q\,\eta\left(\alpha^{\left(1\right)}\right),
    \qquad
    \eta\left(\alpha^{\left(3\right)}\right)
    =
    q^{-i+1}\left[i-1\right]_q\,\eta\left(\alpha^{\left(1\right)}\right),
\]
and the exponents of \(\left(-q\right)\) satisfy
\[
    d_{\alpha^{\left(2\right)}}\left(w\right)-\ell\left(\alpha^{\left(2\right)}\right)
    =
    d_{\alpha^{\left(1\right)}}\left(w\right)-\ell\left(\alpha^{\left(1\right)}\right)+1,
\]
while
\[
    d_{\alpha^{\left(3\right)}}\left(w\right)-\ell\left(\alpha^{\left(3\right)}\right)
    =
    d_{\alpha^{\left(1\right)}}\left(w\right)-\ell\left(\alpha^{\left(1\right)}\right).
\]
Thus the sum of the three local coefficients is
\[
\begin{aligned}
    1
    +\left(-q\right)q^{-i}\left[i\right]_q
    +q^{-i+1}\left[i-1\right]_q
    &=
    1-q^{-i+1}\left(\left[i-1\right]_q+q^{i-1}\right)
      +q^{-i+1}\left[i-1\right]_q=0.
\end{aligned}
\]
Therefore every word with a descent has total coefficient zero.

It remains to consider the increasing word
\(\left(a_1,a_2,\ldots,a_r\right)\). This word is \(\alpha\)-\(V\)-shaped for every
\(\alpha\models r\), and in every block the valley occurs at the first position.
Hence \(d_\alpha\left(w\right)=\ell\left(\alpha\right)\). Its total coefficient is therefore
\[
    q^{\binom r2}
    \sum_{\alpha\models r}\eta\left(\alpha\right)
    =
    q^{\binom r2}
    \sum_{D\subseteq \left[r-1\right]}
    \prod_{i\in \left[r-1\right]\setminus D}q^{-i}\left[i\right]_q .
\]
The sum over $D\subseteq[r-1]$ factors as
\[
    \prod_{i=1}^{r-1}\left(1+q^{-i}\left[i\right]_q\right)
    =
    \prod_{i=1}^{r-1}q^{-i}\left[i+1\right]_q
    =
    q^{-\binom r2}\left[r\right]_q!.
\]
Thus the coefficient of
\(\left\langle a_1,a_2,\ldots,a_r\right\rangle\) is \(\left[r\right]_q!\), and all other
coefficients vanish. 
\end{proof}

\begin{example}
    Let $A=\left\{1,2\right\}$. Then the left-hand side of Proposition~\ref{prop:technical-bracket} is
    \begin{equation*}
        S\left(1,1\right)\sum_{W\in \mathcal{B}\left(A;\left(1,1\right)\right)}\left\langle W \right\rangle+S\left(2\right)\sum_{W\in \mathcal{B}\left(A;\left(2\right)\right)}\left\langle W \right\rangle.
    \end{equation*}
Substituting $S\left(1,1\right)=q$ and $S\left(2\right)=1-q$, we obtain
\begin{align*}
    &q\sum_{W\in\mathcal B\left(A;\left(1,1\right)\right)}\left\langle W\right\rangle
    +\left(1-q\right)\sum_{W\in\mathcal B\left(A;\left(2\right)\right)}\left\langle W\right\rangle\\
    &\qquad=
    q\left\langle 1,2\right\rangle+q\left\langle 2,1\right\rangle
    +\left(1-q\right)\left\langle \left[2,1\right]\right\rangle\\
    &\qquad=
    q\left\langle 1,2\right\rangle+q\left\langle 2,1\right\rangle
    +\left\langle 1,2\right\rangle-q\left\langle 2,1\right\rangle\\
    &\qquad=
    \left(1+q\right)\left\langle 1,2\right\rangle,
\end{align*}
where the second equality follows from
Lemma~\ref{lem:q-swapping}.
\end{example}

\begin{proof}[Proof of \eqref{eq: dual llt}]
  First write
   \begin{equation*}
       \widetilde{\Xi}_{\mmb}\left[X;q\right]=\sum_{\alpha\models k}\sum_{W\in \mathcal{B}\left(\alpha\right)}\left\langle W \right\rangle M_{\alpha}.
   \end{equation*}
    Fix a partition
    $\lambda=(\lambda_1,\ldots,\lambda_\ell)\vdash k$.
    By Corollary~\ref{cor: m plethysm}, we have
   \begin{equation*}
       \left[m_{\lambda}\right] \left(\widetilde{\Xi}_{\mmb}\left[\frac{X}{1-q};q\right]\right)=\frac{1}{\left(1-q\right)^{\left|\lambda\right|}\left[\lambda\right]_{q}!}\left(\sum_{\substack{\alpha^{\left(1\right)},\alpha^{\left(2\right)},\dots,\alpha^{\left(\ell\right)}\\\alpha^{\left(i\right)}\models \lambda_i}}\prod_{i=1}^{\ell}S\left(\alpha^{\left(i\right)}\right)\sum_{W\in \mathcal{B}\left(\alpha^{\left(1\right)}\alpha^{\left(2\right)}\cdots\alpha^{\left(\ell\right)}\right)}\left\langle W \right\rangle\right).
   \end{equation*}
For each $i$, the bracketed blocks specified by $\alpha^{\left(i\right)}$
occupy a consecutive segment of $W$ of length $\lambda_i$.
We group the terms in the sum over
$\alpha^{\left(1\right)},\ldots,\alpha^{\left(\ell\right)}$ according to the sets
of letters in these segments. Thus let
\[
    \left[k\right]=A_1\sqcup\cdots\sqcup A_\ell,
    \qquad \left|A_i\right|=\lambda_i,
\]
and write
\[
    A_i=\left\{a^{\left(i\right)}_1<\cdots<a^{\left(i\right)}_{\lambda_i}\right\}.
\]
Applying Proposition~\ref{prop:technical-bracket}
successively to these segments, we find that the contribution of this
ordered set partition to the sum is
\[
    \left[\lambda\right]_q!\,
    \left\langle
        a^{\left(1\right)}_1,\ldots,a^{\left(1\right)}_{\lambda_1},
        a^{\left(2\right)}_1,\ldots,a^{\left(2\right)}_{\lambda_2},
        \ldots,
        a^{\left(\ell\right)}_1,\ldots,a^{\left(\ell\right)}_{\lambda_\ell}
    \right\rangle.
\]
Hence we conclude
\begin{equation}\label{eq:dual-rllt-sorted-coefficient}
    \left(1-q\right)^k
    \left[m_\lambda\left[X\right]\right]
    \widetilde{\Xi}_{\mm^\bullet}\!\left[\frac{X}{1-q};q\right]  =
    \sum_{\substack{
        \left[k\right]=A_1\sqcup\cdots\sqcup A_\ell\\
        \left|A_i\right|=\lambda_i
    }}
    \left\langle
        A_1^{\uparrow},A_2^{\uparrow},\ldots,A_\ell^{\uparrow}
    \right\rangle ,
\end{equation}
where \(A_i^{\uparrow}\) denotes the increasing listing of \(A_i\). Note that the right-hand side of \eqref{eq:dual-rllt-sorted-coefficient} is precisely
\begin{equation*}
    \sum_{L\in \overline{\sLD}\left(\mmb;\lambda\right)}q^{\asc\left(L\right)}.
\end{equation*}
We may replace $\lambda$ with any composition $\beta$ such that $\beta\sim\lambda$. Thus, the right-hand side of \eqref{eq:dual-rllt-sorted-coefficient} also equals
    \begin{equation*}
       \sum_{L\in \overline{\sLD}\left(\mmb;\beta\right)}q^{\asc\left(L\right)}.
   \end{equation*}

    We finally conclude that 
   \begin{equation*}
       \left(1-q\right)^k \widetilde{\Xi}_{\mmb}\left[\frac{X}{1-q};q\right]=\sum_{\lambda \vdash k}\sum_{L\in \overline{\sLD}\left(\mm^{\bullet};\lambda\right)}q^{\asc\left(L\right)}m_{\lambda}=\sum_{\beta \models k}\sum_{L\in \overline{\sLD}\left(\mm^{\bullet};\beta\right)}q^{\asc\left(L\right)}M_{\beta}.
   \end{equation*}
\end{proof}

\begin{proof}[Proof of \eqref{eq: llt}]
For $L=\left(L_1,\ldots,L_k\right)\in\LD\left(\mmb\right)$, define
\[
    D_1\left(L\right)\coloneqq\left\{i\in\left[k-1\right]:L_i\left(1\right)>L_{i+1}\left(1\right)\right\}.
\]
Since $\left(L_1\left(1\right),\ldots,L_k\left(1\right)\right)$ is a permutation of $\left[k\right]$,
for every composition $\beta\models k$ we have
\[
    L\in\overline{\sLD}\left(\mmb;\beta\right)
    \quad\Longleftrightarrow\quad
    D_1\left(L\right)\subseteq\Set\left(\beta\right).
\]
Therefore, the monomial quasisymmetric expansion obtained
in the proof of \eqref{eq: dual llt} can be rewritten as
\[
    \left(1-q\right)^k
    \widetilde{\Xi}_{\mmb}\!\left[\frac{X}{1-q};q\right]
    =
    \sum_{L\in\LD\left(\mmb\right)}
        q^{\asc\left(L\right)}F_{D_1\left(L\right)}\left[X\right].
\]

On symmetric functions, $\omega$ commutes with the
plethystic substitution $X\mapsto X/\left(1-q\right)$, since both
operators are diagonal on the power-sum basis.
By Lemma~\ref{lem:omega-rational-chromatic}, we also have
$\omega\widetilde{\Xi}_{\mmb}=\Xi_{\mmb}$.
Applying $\omega$ to the expansion in $F_{D_1\left(L\right)}\left[X\right]$
and using \eqref{eq: involution on F}, we obtain
\[
    \left(1-q\right)^k
    \Xi_{\mmb}\!\left[\frac{X}{1-q};q\right]
    =
    \sum_{L\in\LD\left(\mmb\right)}
        q^{\asc\left(L\right)}F_{D_1\left(L\right)^c}\left[X\right],
\]
where $D_1\left(L\right)^c=\left[k-1\right]\setminus D_1\left(L\right)$.

For a composition $\beta\models k$, the condition
$D_1\left(L\right)^c\subseteq\Set\left(\beta\right)$ is equivalent to
$L\in\sLD\left(\mmb;\beta\right)$.
Expanding the fundamental quasisymmetric functions
therefore gives
\[
    \left(1-q\right)^k
    \Xi_{\mmb}\!\left[\frac{X}{1-q};q\right]
    =
    \sum_{\beta\models k}
    \sum_{L\in\sLD\left(\mmb;\beta\right)}
        q^{\asc\left(L\right)}M_\beta\left[X\right].
\]
Since the left-hand side is symmetric, the coefficients of
$M_\beta\left[X\right]$ agree for rearrangements of $\beta$.
Grouping these terms into monomial symmetric functions
proves \eqref{eq: llt}.
\end{proof}

\subsection{Proof of Theorem~\ref{thm: main B}}

We first define the terms used in Theorem~\ref{thm: main B}.
\begin{definition}
Let $(\pi^\bullet,b)$ belong to either
$\cPF_{m,n}^k\times\mathbb Z_{>0}^k$ or
$\overline{\cPF}_{m,n}^k\times\mathbb Z_{>0}^k$.
We call $(\pi^\bullet,b)$ \emph{sorted} if
\[
    b_1\le b_2\le\cdots\le b_k
\]
and, whenever $b_i=b_{i+1}$, one has
\[
    \pi^{(i)}_1\le\pi^{(i+1)}_1
\]
with respect to the total order on labeled north steps.

For either type of tuple, define
\[
    \aut_q\left(\pi^\bullet,b\right)\coloneqq\left[\alpha\right]_q!,
\]
where \(\alpha\) records the multiplicities of the elements in the multiset
\[
    \left\{\left(\pi^{\left(i\right)}_1,b_i\right):1\le i\le k\right\}.
\]
\end{definition}

\begin{thm}\label{thm: cm general}
Let $m,n$ be coprime positive integers, and let $k\ge1$. Then
\begin{align}
    e_k[-X^{m,n}Y]\cdot1
    &=
    \sum_{\substack{
        \left(\pi^\bullet,b\right)\in
        \cPF_{m,n}^k\times\mathbb Z_{>0}^k\\
        \left(\pi^\bullet,b\right)\text{ sorted}
    }}
    \frac{
        t^{\area\left(\pi^\bullet\right)}
        q^{\stat\left(\pi^\bullet\right)}
    }{
        \left(1-q\right)^k\aut_q\left(\pi^\bullet,b\right)
    }
    X_{\pi^\bullet}Y_b ,
    \label{eq:cm-general-strict}\\
    \omega_X\left(e_k[-X^{m,n}Y]\cdot1\right)
    &=
    \sum_{\substack{
        \left(\pi^\bullet,b\right)\in
        \overline{\cPF}_{m,n}^k\times\mathbb Z_{>0}^k\\
        \left(\pi^\bullet,b\right)\text{ sorted}
    }}
    \frac{
        t^{\area\left(\pi^\bullet\right)}
        q^{\overline{\stat}\left(\pi^\bullet\right)}
    }{
        \left(1-q\right)^k\aut_q\left(\pi^\bullet,b\right)
    }
    X_{\pi^\bullet}Y_b.
    \label{eq:cm-general-weak}
\end{align}
Here, we write
$Y_b\coloneqq y_{b_1}\cdots y_{b_k}$ for $b=(b_1,\ldots,b_k)\in\mathbb Z_{>0}^k$.

\end{thm}

\begin{proof}
We prove \eqref{eq:cm-general-weak}. Equation~\eqref{eq:cm-general-strict} can be proved in a similar manner. 

By Theorem~\ref{thm: main} and Theorem~\ref{thm: rLLt monomial}, we have
\begin{equation}\label{eq:weak-after-rllt}
\begin{aligned}
    \omega_X\left(e_k[-X^{m,n}Y]\cdot1\right)
    =
    \frac{1}{\left(1-q\right)^k}
    \sum_{\alpha\models kn}
    M_\alpha\left[X\right]
    \sum_{\substack{
        \tau^\bullet\in
        \overline{\cPF}_{m,n}^k\left(\alpha\right)\\
        \tau^\bullet\text{ dominant}
    }}
    t^{\area\left(\tau^\bullet\right)}
    q^{\overline{\stat}\left(\tau^\bullet\right)}
    \frac{
        \sum_{\beta\models k}
        \sum_{L\in\sLD\left(\mmb\left(\tau^\bullet\right);\beta\right)}
            q^{\asc\left(L\right)}M_\beta\left[Y\right]
    }{
        \left[\alpha'\left(\tau^\bullet\right)\right]_q!
    } .
\end{aligned}
\end{equation}

For a composition \(\beta=\left(\beta_1,\ldots,\beta_\ell\right)\models k\), let
\(\mathbb Z_{>0}^k\left(\beta\right)\) be the set of vectors
\(b=\left(b_1,\ldots,b_k\right)\) satisfying
\[
    b_1=\cdots=b_{\beta_1}
    <
    b_{\beta_1+1}=\cdots=b_{\beta_1+\beta_2}
    <\cdots .
\]
Then $M_\beta\left[Y\right]
    =
    \sum_{b\in\mathbb Z_{>0}^k\left(\beta\right)}Y_b$.

Fix a dominant \(\tau^\bullet\), a composition \(\beta\models k\), a loop
decomposition \(L=\left(L_1,\ldots,L_k\right)\in\sLD\left(\mmb\left(\tau^\bullet\right);\beta\right)\), and
\(b\in\mathbb Z_{>0}^k\left(\beta\right)\). Define
\(\rho^\bullet=\rho^\bullet\left(\tau^\bullet,L\right)\) by
\[
    \rho^{\left(i\right)}_\ell\coloneqq\tau^{\left(k-L_i\left(\ell\right)+1\right)}_\ell
    \qquad\left(1\le i\le k,\ 1\le\ell\le n\right).
\]
Since $b_i=b_{i+1}$ forces $L_i(1)>L_{i+1}(1)$,
dominance of $\tau^\bullet$ and the reversed indexing in the
definition of $\rho^\bullet$ give
$\rho^{(i)}_1\le\rho^{(i+1)}_1$; hence
$(\rho^\bullet,b)$ is sorted.

Now we compute the coefficient with which a fixed sorted pair
\(\left(\rho^\bullet,b\right)\) occurs. In other words, we compute
\begin{equation*}
\sum_{\substack{L\in\sLD\left(\mmb\left(\tau^\bullet\right);\beta\right)\\\rho^\bullet=\rho^\bullet\left(\tau^\bullet,L\right)}} q^{\asc\left(L\right)}.
\end{equation*}

Let $\gamma\left(i,\ell\right)\coloneqq\tau^{\left(k-i+1\right)}_\ell
$ be the labeled north step attached to the vertex \(\left(i,\ell\right)\) of
\(G\left(\mmb\left(\tau^\bullet\right)\right)\). Partition the vertex set of
\(G\left(\mmb\left(\tau^\bullet\right)\right)\) into equivalence classes \(I\), where two vertices are
equivalent if they have the same image under \(\gamma\). If \(r_I=\left|I\right|\), then
\[
    \left[\alpha'\left(\tau^\bullet\right)\right]_q!
    =
    \prod_I \left[r_I\right]_q!.
\]

For each class $I$, list its vertices in increasing order of their first
coordinates, and let $d_1,\ldots,d_{r_I}$ be the indices of the loops
containing these vertices. As in Proposition~\ref{prop: L to pib}, let
$w^I\in\mathfrak S_{r_I}$ be obtained by standardizing
$(d_{r_I},\ldots,d_1)$. The ascent statistic decomposes as
\begin{equation}\label{eq:asc-pos-inv-local}
    \asc\left(L\right)
    =
    \overline{\pos}\left(\rho^\bullet\right)+
    \sum_I \inv\left(w^I\right).
\end{equation}
This is the same local correspondence as in
Proposition~\ref{prop: L to pib}.

Since \(\tau^\bullet\) is dominant, \(\overline{\pos}(\tau^\bullet)=0\). Since \(\rho^\bullet\) and \(\tau^\bullet\) have the same skeleton, Lemma~\ref{lem: pos crossing} gives
\[
    \overline{\stat}\left(\rho^\bullet\right)
    =
    \overline{\stat}\left(\tau^\bullet\right)
    +\overline{\pos}\left(\rho^\bullet\right).
\]
Combining this with \eqref{eq:asc-pos-inv-local}, we obtain
\begin{equation}\label{eq:q-weight-transfer}
    q^{\overline{\stat}\left(\tau^\bullet\right)}
    q^{\asc\left(L\right)}
    =
    q^{\overline{\stat}\left(\rho^\bullet\right)}
    \prod_I q^{\inv\left(w^I\right)}.
\end{equation}

Now fix the sorted pair \(\left(\rho^\bullet,b\right)\). For an equivalence class
\(I\) not lying at height \(1\), all permutations \(w^I\in\mathfrak S_{r_I}\)
occur, and their inversion generating function is
\[
    \sum_{w^I\in\mathfrak S_{r_I}}q^{\inv\left(w^I\right)}
    =
    \left[r_I\right]_q!.
\]
Thus these factors cancel the corresponding factors in
\(\left[\alpha'\left(\tau^\bullet\right)\right]_q!\).

It remains to consider the equivalence classes $I$ lying at height $1$.
Let $\eta_I$ be the common labeled north-step type of the vertices in $I$,
and put
\[
    J_I=\{i\in[k]:\rho^{(i)}_1=\eta_I\}.
\]
List the distinct values of $b_i$, for $i\in J_I$, in increasing order as
$c_1<\cdots<c_s$, and define
\[
    \nu_a=\#\{i\in J_I:b_i=c_a\},
    \qquad
    N_a=\nu_1+\cdots+\nu_a,
    \qquad N_0=0.
\]
Thus $\nu(I;\rho^\bullet,b)=(\nu_1,\ldots,\nu_s)$ records the
multiplicities in increasing order of the $b$-values.

The values of $w^I$ are the ranks of the loop indices in $J_I$,
whereas its positions follow the vertices of $I$ in decreasing order.
Since the loop vertices decrease within each equal-$b$ block,
the allowed permutations $w^I$ are precisely those satisfying
\[
    (w^I)^{-1}(N_{a-1}+1)
    <\cdots<
    (w^I)^{-1}(N_a)
    \qquad(1\le a\le s).
\]
These are the minimal-length representatives for the quotient
\[
    \left(
        \mathfrak S_{\nu_1}\times\cdots\times\mathfrak S_{\nu_s}
    \right)
    \backslash\mathfrak S_{r_I}.
\]
Inversion identifies these representatives with the minimal-length
representatives for
$\mathfrak S_{r_I}/(\mathfrak S_{\nu_1}\times\cdots\times
\mathfrak S_{\nu_s})$ and preserves the inversion number. Therefore
\[
    \sum_{\text{allowed }w^I}q^{\inv\left(w^I\right)}
    =
    \frac{\left[r_I\right]_q!}
    {\left[\nu_1\right]_q!\cdots\left[\nu_s\right]_q!}.
\]
After division by the factor \(\left[r_I\right]_q!\) coming from
\(\left[\alpha'\left(\tau^\bullet\right)\right]_q!\), the contribution of this class is
\[
    \frac{1}{\left[\nu_1\right]_q!\cdots\left[\nu_s\right]_q!}.
\]
Multiplying over all equivalence classes at height \(1\), we get exactly
\[
    \frac{1}{\aut_q\left(\rho^\bullet,b\right)}.
\]

Using \eqref{eq:q-weight-transfer} together with
\(\area\left(\rho^\bullet\right)=\area\left(\tau^\bullet\right)\)
and \(X_{\rho^\bullet}=X_{\tau^\bullet}\), we find that the contribution of all data mapping to a fixed sorted pair
\(\left(\rho^\bullet,b\right)\) is
\[
    \frac{
        t^{\area\left(\rho^\bullet\right)}
        q^{\overline{\stat}\left(\rho^\bullet\right)}
    }{
        \left(1-q\right)^k\aut_q\left(\rho^\bullet,b\right)
    }
    X_{\rho^\bullet}Y_b.
\]
Summing over all sorted pairs gives
\[
    \omega_X\left(e_k[-X^{m,n}Y]\cdot1\right)
    =
    \sum_{\substack{
        \left(\rho^\bullet,b\right)\in
        \overline{\cPF}_{m,n}^k\times\mathbb Z_{>0}^k\\
        \left(\rho^\bullet,b\right)\text{ sorted}
    }}
    \frac{
        t^{\area\left(\rho^\bullet\right)}
        q^{\overline{\stat}\left(\rho^\bullet\right)}
    }{
        \left(1-q\right)^k\aut_q\left(\rho^\bullet,b\right)
    }
    X_{\rho^\bullet}Y_b,
\]
which is \eqref{eq:cm-general-weak}.
\end{proof}

\begin{rmk}
When $n=1$, equation~\eqref{eq:eha-nabla-identity} gives
\[
    e_k[-X^{m,1}Y]\cdot1
    =\nabla^m e_k\!\left[\frac{XY}{M}\right].
\]
Thus Theorem~\ref{thm: cm general} gives a cyclic
parking-function expression for the same symmetric function
as in~\cite[Theorem~A]{CM25}. Note that \cite[Proposition 3.2]{CM25} expresses this symmetric function in terms of the unicellular LLT polynomials. 
\end{rmk}

\appendix

\section{Affine permutations and cyclic parking functions}\label{Sec: Appendix}

\subsection{(Extended) affine permutations}
In this subsection, we define affine permutations and the necessary terminology. We then recall the area and \(\dinv\) statistics on affine permutations from \cite{CM25}.

    An (extended) affine permutation $w$ of size $n$ is a bijection $\mathbb{Z}\rightarrow \mathbb{Z}$ satisfying
    \begin{equation*}
        w\left(i+n\right)=w\left(i\right)+n
    \end{equation*}
    for all $i\in \mathbb{Z}$. Such a $w$ is often written in window notation as
    \begin{equation*}
       \left[w\left(1\right),w\left(2\right),\dots,w\left(n\right)\right],
    \end{equation*}
    which determines $w$ uniquely. We denote the set of affine permutations of size $n$ by $W_n$. We emphasize that we \emph{do not} impose the condition $\sum_{i=1}^{n}w\left(i\right)=\frac{n\left(n+1\right)}{2}$, which is often assumed in other contexts, since we are working with \emph{extended} affine permutations. We say that an affine permutation is \emph{positive} if $w\left(x\right)>0$ for all $x>0$.

    For a positive integer $m$, an affine permutation $w$ is called \emph{$m$-restricted} if $w^{-1}\left(x+m\right)>w^{-1}\left(x\right)$ for all $x\in \mathbb{Z}$. We denote by $S_{m,n}^{+}$ the set of $m$-restricted positive affine permutations of size $n$.

    We define
    \[
        \area\left(w\right)\coloneqq\sum_{i=1}^{n}\frac{w\left(i\right)-i}{n}.
    \]
    For $a,b\in\mathbb Z$ with $a\not\equiv b\pmod n$, let $t_{a,b}$ be the affine transposition defined by
    \[
    t_{a,b}\left(x\right)=
    \begin{cases}
    x+\left(b-a\right), & x\equiv a \pmod n,\\
    x+\left(a-b\right), & x\equiv b \pmod n,\\
    x,       & x\not\equiv a \pmod n \text{ and } x\not\equiv b \pmod n.
    \end{cases}
    \]
    We define
    \begin{equation*}
        E_{m}\left(w\right)=\left\{t_{a,b}:t_{a,b}w < w, \left|a-b\right|<m\right\},
    \end{equation*}
     where $<$ between affine permutations is the usual affine Bruhat order. Following~\cite[Section~5.1, equation~(39)]{CM25},
    we use this set of affine Bruhat reflections to define
    the statistic $\dinv_m$ below. Finally, define 
    \begin{equation*}
        \dinv_m\left(w\right)=\frac{mn-m-n+\gcd\left(m,n\right)}{2}-\left|E_m\left(w\right)\right|.
    \end{equation*}

Let \(m,n\) be coprime positive integers, and let \(k\ge 1\).
Following Carlsson--Mellit \cite[Equation~(1)]{CM25}, set
\[
\mathcal{H}_{km,kn}
=
\frac{1}{\left(1-q\right)^k}
\sum_{w\in S^{+}_{km,kn}}
t^{\area\left(w\right)}q^{\dinv_{km}\left(w\right)}.
\]
Carlsson--Mellit \cite[Corollary~5.7]{CM25} proved the $n=1$
case of the coefficient identity
\[
    \mathcal H_{km,kn}
    =
    \left[m_{1^{kn}}\left[X\right]m_{1^k}\left[Y\right]\right]
    \left(e_k\!\left[-X^{m,n}Y\right]\cdot 1\right).
\]
In the following subsection, we establish this identity for arbitrary
coprime positive integers $m,n$ by combining a statistic-preserving
bijection with the generating-function formula in
\cite[Theorem~1.1]{KO26}; see Corollary~\ref{cor: Hilbert series}.

\subsection{A statistic-preserving bijection}\label{subsec: affine permutation vs cPF}

The relation between affine permutations and labeled rational paths goes
back to Gorsky--Mazin--Vazirani \cite{GMV16}; invariant-set and multi-path
models in the noncoprime case are developed in \cite{GMV17,GMV20}.
The normalized affine permutations considered in \cite{GMV16}
do not comprise the full set of positive extended affine permutations
used here; the rational parking functions in that work correspond
to the normalized objects. We give the explicit cyclic-tuple correspondence, including its
inverse and both weight identities, for completeness.
Together with \cite[Theorem~1.1]{KO26}, it will yield
Corollary~\ref{cor: Hilbert series}.

Given $w\in S^+_{km,kn}$ and $1\le i\le kn$, place a north step of content
$w\left(i\right)-1$ in the strip $0\le y<n$ and label it by $i$.
The content determines a unique sheet, height, and horizontal coordinate:
its residue modulo $k$ determines the sheet, and then coprimality of $m,n$
determines the height modulo $n$.
Connect successive north steps in each sheet by east steps, with the last
one connected to the $\left(m,n\right)$-translate of the first. Denote the resulting
tuple by $\phi\left(w\right)$.

\begin{lem}\label{lem: affine cpf bijection}
The map $\phi:S^+_{km,kn}\longrightarrow\cPF_{m,n}^{k,\mathrm{std}}$ is a
bijection, and
\[
 \area\left(w\right)=\area\left(\phi\left(w\right)\right),\qquad
 \dinv_{km}\left(w\right)=\stat\left(\phi\left(w\right)\right).
\]
\end{lem}
\begin{proof}
We first show that $\phi\left(w\right)$ is a $k$-tuple of cyclic parking
functions. For a north step $N\left(i\right)$ labeled $i$, let
$1\leq j\leq kn$ be the unique index such that
$w\left(j\right)\equiv w\left(i\right)+km\pmod{kn}$.
Then $N\left(j\right)$ is the north step in the same sheet at the next
height, interpreted cyclically.
By periodicity, there is a unique integer $a$ such that
$w\left(j+kna\right)=w\left(i\right)+km$.
The $km$-restricted condition gives $j+kna>i$, so $a\geq0$.
If $a=0$, then $N\left(i\right)$ and $N\left(j\right)$ are vertically adjacent
with strictly increasing labels, since $j>i$.
If $a>0$, then the appropriate cyclic translate of $N\left(j\right)$
starts strictly to the right of the endpoint of $N\left(i\right)$.
Thus these two labeled north steps satisfy the cyclic
parking function condition.

We describe the inverse map $\phi^{-1}$.
Given $\pi^\bullet\in\cPF_{m,n}^{k,\std}$, let $N\left(i\right)$ be
the north step labeled $i$ for $1\leq i\leq kn$, and define
$w=\phi^{-1}\left(\pi^\bullet\right)$ by
\[
    w\left(i\right)=\cont\left(N\left(i\right)\right)+1
    \qquad \left(1\leq i\leq kn\right),
\]
extending by $w\left(i+kna\right)=w\left(i\right)+kna$ for all $a\in\mathbb Z$. The north-step contents are pairwise incongruent modulo $kn$:
their residues modulo $k$ distinguish the sheets, and coprimality of
$m,n$ distinguishes the heights within each sheet.
Thus $w$ is an extended affine permutation. Moreover, all north-step
contents are nonnegative because the paths stay weakly above
$my=nx$, so $w$ is positive. 

We show that $w$ is $km$-restricted.
Suppose, for a contradiction, that
$w^{-1}\left(c+km\right)<w^{-1}\left(c\right)$ for some $c\in\mathbb Z$.
Choose $1\leq i\leq kn$ such that
$c'\coloneqq w\left(i\right)\equiv c\pmod{kn}$.
By periodicity, we also have
$w^{-1}\left(c'+km\right)<w^{-1}\left(c'\right)=i$.
Set $i'\coloneqq w^{-1}\left(c'+km\right)$.
If $1\leq i'\leq kn$, then $N\left(i\right)$ and $N\left(i'\right)$ are
vertically adjacent, interpreted cyclically, but their
labels satisfy $i'<i$, contradicting the cyclic labeling
condition.
If $i'<1$, let $1\leq j\leq kn$ be the representative
of $i'$ modulo $kn$.
Then the appropriate cyclic translate of $N\left(j\right)$ starts
strictly to the left of the endpoint of $N\left(i\right)$, again
a contradiction.
Therefore $w$ is $km$-restricted.
The two constructions are inverse to each other, so
$\phi$ is a bijection.

The equality $\area\left(w\right)=\area\left(\phi\left(w\right)\right)$ follows from the
fact that the number of full cells to the right of $N\left(i\right)$
and to the left of the line $my=nx$ is
\[
    \frac{w\left(i\right)-1-b_i}{kn},
\]
where $b_i$ is the least nonnegative residue of $w\left(i\right)-1$
modulo $kn$.
Since $\left\{b_1,\ldots,b_{kn}\right\}=\left\{0,\ldots,kn-1\right\}$, we have
\[
    \area\left(\phi\left(w\right)\right)
    =
    \sum_{i=1}^{kn}\frac{w\left(i\right)-1-b_i}{kn}
    =
    \sum_{i=1}^{kn}\frac{w\left(i\right)-i}{kn}
    =
    \area\left(w\right).
\]

We now show that $\dinv_{km}\left(w\right)=\stat\left(\phi\left(w\right)\right)$.
Let $S\coloneqq\left\{w\left(i\right)-1:1\leq i\leq kn\right\}$ be the set
of north-step contents.
Every affine transposition in $E_{km}\left(w\right)$ has a unique
representative of the form
\[
    t_{s+1,h+1},
    \qquad
    s\in S,\qquad
s<h<s+km,\qquad
h\not\equiv s\pmod{kn}.
\]
By the affine Bruhat reflection criterion,
\[
    t_{s+1,h+1}w<w
    \quad\Longleftrightarrow\quad
    w^{-1}\left(s+1\right)>w^{-1}\left(h+1\right).
\]
We consider three cases according to $h$.

First, suppose $h$ is the content of a cell below
$P_{\phi\left(w\right)}$. Let $u$ be the north-step content in the same row and sheet.
Then $h=u-kn a$ for some integer $a\geq1$.
If $j$ is the label of that north step, periodicity gives
\[
w^{-1}(h+1)=j-kn a\leq0<w^{-1}(s+1).
\]
Thus the affine Bruhat reflection criterion is satisfied.
There is a unique nonnegative integer $\ell$ such that
$h+km\ell$ is the content of an east step of $P_{\phi\left(w\right)}$.
The pair consisting of the north step of content $s$
and this east step is counted by $\pdinv\left(P_{\phi\left(w\right)}\right)$.
Conversely, given a pair $\left(N,E\right)$ counted by
$\pdinv\left(P_{\phi\left(w\right)}\right)$, choose the unique nonnegative
integer $\ell$ such that
\[
    \cont\left(N\right)<\cont\left(E\right)-km\ell<\cont\left(N\right)+km.
\]
Then
$t_{\cont\left(N\right)+1,\cont\left(E\right)-km\ell+1}\in E_{km}\left(w\right)$.

Next, suppose $h\in S$.
Then $s$ and $h$ are both north-step contents, and
$s<h<s+km$ is exactly the attacking-pair condition.
Moreover, $w^{-1}\left(s+1\right)>w^{-1}\left(h+1\right)$ means that the
labels of this attacking pair form an inversion,
so the pair is counted by $\ldinv\left(\phi\left(w\right)\right)$.
Conversely, every pair counted by $\ldinv\left(\phi\left(w\right)\right)$
gives an element of $E_{km}\left(w\right)$ in this way.

Finally, suppose $h$ is the content of a cell above
$P_{\phi\left(w\right)}$ and is not a north-step content.
Then there is a positive integer $\ell$ such that
$h-kn\ell\in S$.
If $u$ is the label of the north step with content
$h-kn\ell$, then
\[
    w^{-1}\left(h+1\right)=u+kn\ell.
\]
Therefore
\[
    w^{-1}\left(h+1\right)>kn\ell\geq kn\geq w^{-1}\left(s+1\right),
\]
so no element of $E_{km}\left(w\right)$ occurs in this case.

Consequently,
\[
    \left|E_{km}\left(w\right)\right|
    =
    \pdinv\left(P_{\phi\left(w\right)}\right)+\ldinv\left(\phi\left(w\right)\right),
\]
which proves the desired identity.
\end{proof}

\begin{example}\label{ex:affine-cyclic-pf-k2-m3-n2}
Let $m=3,n=2,k=2$ and consider the positive affine permutation $w=\left[1,4,7,6\right]$ of size $4$,
which is $6$-restricted. Then $\phi\left(w\right)\in \cPF_{3,2}^{2,\std}$ is obtained by placing north steps
with contents $0,3,6,5$ and labels $1,2,3,4$, respectively. Figure~\ref{fig:affine-cyclic-pf-k2-m3-n2} shows $\phi\left(w\right)=\left(\pi^{\left(1\right)},\pi^{\left(2\right)}\right)$. 

\begin{figure}[ht]
\centering

$\displaystyle \pi^{\left(1\right)}=$
\begin{tikzpicture}[scale=0.85, baseline=(current bounding box.center)]
    \foreach \x in {-1,...,2} {
        \foreach \y in {0,...,1} {
            \pgfmathtruncatemacro{\val}{6*(\y)-4*(\x+1)}

            \ifnum\val<0
                \def\txtstyle{}
            \else
                \def\txtstyle{\bfseries}
            \fi

            \def\txtcolor{black!60}

            \ifnum\y=0
                \ifnum\x=-1
                    \def\txtcolor{green!60!black}
                \fi
            \fi
            \ifnum\y=1
                \ifnum\x=-1
                    \def\txtcolor{green!60!black}
                \fi
            \fi

            \ifnum\y=1
                \ifnum\x=0
                    \def\txtcolor{green!60!black}
                \fi
            \fi
            \ifnum\y=1
                \ifnum\x=1
                    \def\txtcolor{green!60!black}
                \fi
            \fi
            \ifnum\y=1
                \ifnum\x=2
                    \def\txtcolor{green!60!black}
                \fi
            \fi

            \node[\txtcolor, font=\txtstyle, scale=0.8]
                at (\x+0.5,\y+0.5) {\val};
        }
    }

    \draw[help lines, gray!30, thin] (-1,0) grid (3,2);
    \draw[black, thin] (0,0) -- (3,2);

    \draw[red, line width=2pt, line cap=round, line join=round]
        (0,0) -- (0,1) -- (0,2)
        -- (1,2) -- (2,2) -- (3,2);

    \node[fill=white, inner sep=1pt, text=red, font=\scriptsize]
        at (0,0.5) {$1$};
    \node[fill=white, inner sep=1pt, text=red, font=\scriptsize]
        at (0,1.5) {$3$};
\end{tikzpicture}
\qquad
$\displaystyle \pi^{\left(2\right)}=$
\begin{tikzpicture}[scale=0.85, baseline=(current bounding box.center)]
    \foreach \x in {-2,...,2} {
        \foreach \y in {0,...,1} {
            \pgfmathtruncatemacro{\val}{1+6*(\y)-4*(\x+1)}

            \ifnum\val<0
                \def\txtstyle{}
            \else
                \def\txtstyle{\bfseries}
            \fi

            \def\txtcolor{black!60}

            \ifnum\y=0
                \ifnum\x=-2
                    \def\txtcolor{green!60!black}
                \fi
            \fi
            \ifnum\y=1
                \ifnum\x=0
                    \def\txtcolor{green!60!black}
                \fi
            \fi

            \ifnum\y=0
                \ifnum\x=-1
                    \def\txtcolor{green!60!black}
                \fi
            \fi
            \ifnum\y=0
                \ifnum\x=0
                    \def\txtcolor{green!60!black}
                \fi
            \fi
            \ifnum\y=1
                \ifnum\x=1
                    \def\txtcolor{green!60!black}
                \fi
            \fi

            \node[\txtcolor, font=\txtstyle, scale=0.8]
                at (\x+0.5,\y+0.5) {\val};
        }
    }

    \draw[help lines, gray!30, thin] (-2,0) grid (3,2);
    \draw[black, thin] (0,0) -- (3,2);

    \draw[blue, line width=2pt, line cap=round, line join=round]
        (-1,0) -- (-1,1) -- (0,1)
        -- (1,1) -- (1,2) -- (2,2);

    \node[fill=white, inner sep=1pt, text=blue, font=\scriptsize]
        at (-1,0.5) {$4$};
    \node[fill=white, inner sep=1pt, text=blue, font=\scriptsize]
        at (1,1.5) {$2$};
\end{tikzpicture}

\caption{The tuple of cyclic parking functions $\phi\left(w\right)\in\cPF_{3,2}^{2,\std}$ for
$w=\left[1,4,7,6\right]$.}
\label{fig:affine-cyclic-pf-k2-m3-n2}
\end{figure}
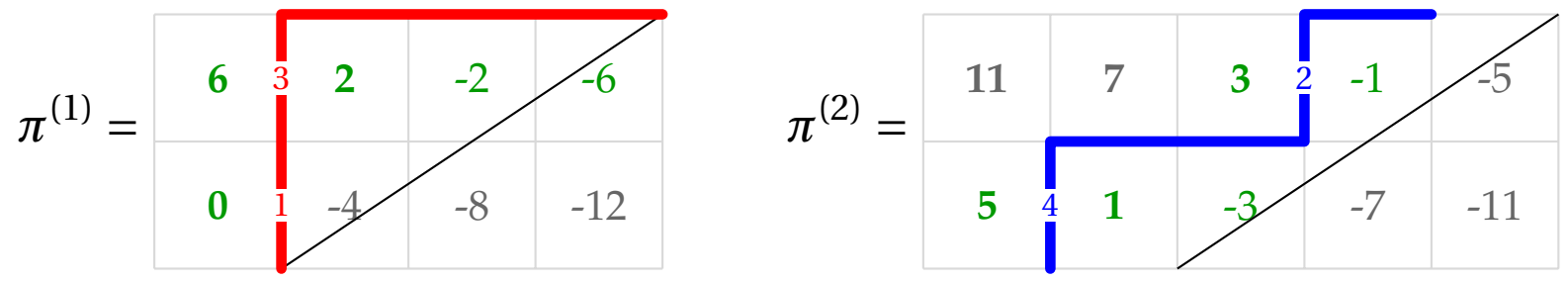
We have $\area\left(w\right)=\frac{1}{4}\sum_{i=1}^{4} \left(w\left(i\right)-i\right)=2=\area(P_{\phi(w)}).$
The set of east-step contents is
$\left\{-6,-3,-2,-1,1,2\right\}$, and the set of north-step contents is
$\left\{0,3,5,6\right\}$. The pairs, listed by their contents, contributing to $\pdinv$ are $\left(0,1\right)$ and $\left(0,2\right)$. Thus, $\pdinv\left(P_{\phi\left(w\right)}\right)=2$.

For the labeled diagonal inversions, the only attacking pair whose labels form an inversion
is the pair of north steps with contents $5$ and $6$, carrying labels $4$ and $3$,
respectively. Hence $\ldinv\left(\phi\left(w\right)\right)=1$.
Finally, $E_{6}\left(w\right)$ is given by $\left\{t_{1,2},t_{1,3},t_{6,7}\right\}$. Thus $\stat\left(\phi\left(w\right)\right)=8-2-1=5$, while
$\dinv_6\left(w\right)=8-\left|E_6\left(w\right)\right|=5$.
\end{example}

\begin{corollary}\label{cor: Hilbert series}
Let $m,n,k$ be positive integers with $\gcd\left(m,n\right)=1$. We have
\begin{equation*}
     \mathcal{H}_{km,kn}=\left[m_{1^{kn}}\left[X\right]m_{1^k}\left[Y\right]\right]\left(e_k[-X^{m,n}Y]\cdot1\right).
\end{equation*}    
\end{corollary}
\begin{proof}
    We have $\left[m_{1^k}\left[Y\right]\right]e_k[-X^{m,n}Y]\cdot1=e_{1^k}\left[-X^{m,n}\right]\cdot 1=\frac{e_{1^k}\left[-MX^{m,n}\right]\cdot 1}{\left(1-q\right)^k\left(1-t\right)^k}$. By \cite[Theorem~1.1]{KO26}, we have 
    \begin{equation*}
        \frac{e_{1^k}\left[-MX^{m,n}\right]\cdot 1}{\left(1-t\right)^k}=\sum_{\pi^{\bullet}\in \cPF_{m,n}^{k}}t^{\area\left(\pi^{\bullet}\right)}q^{\stat\left(\pi^{\bullet}\right)}X_{\pi^{\bullet}}.
    \end{equation*}
Taking the coefficient of $m_{1^{kn}}\left[X\right]$ gives
\[
    \left[m_{1^{kn}}\left[X\right]m_{1^k}\left[Y\right]\right]\left(e_k[-X^{m,n}Y]\cdot1\right)
    =
    \frac{1}{\left(1-q\right)^k}
    \sum_{\pi^\bullet\in\cPF_{m,n}^{k,\mathrm{std}}}
    t^{\area\left(\pi^\bullet\right)}q^{\stat\left(\pi^\bullet\right)}.
\]
By Lemma~\ref{lem: affine cpf bijection}, this equals
\[
    \frac{1}{\left(1-q\right)^k}
    \sum_{w\in S^+_{km,kn}}
    t^{\area\left(w\right)}q^{\dinv_{km}\left(w\right)}
    =
    \mathcal H_{km,kn}.
\]
\end{proof}

\bibliographystyle{alpha}  
\bibliography{main}

\end{document}